\documentclass[a4paper]{amsart}

\usepackage[matrix,arrow,curve,tips]{xy}
            \SelectTips{cm}{}

\newtheorem{dfn}{Definition}[section]
\newtheorem{prop}[dfn]{Proposition}
\newtheorem{theo}[dfn]{Theorem}
\newtheorem{cor}[dfn]{Corollary}
\newtheorem{ex}[dfn]{Example}
\newtheorem{lem}[dfn]{Lemma}

\newcommand{\RR}{\mathbb{R}}
\newcommand{\CC}{\mathbb{C}}
\newcommand{\FF}{\mathbb{F}}

\newcommand{\ra}{\rightarrow}
\newcommand{\lra}{\longrightarrow}
\newcommand{\Rem}{\noindent {\it Remark. }}

\newcommand{\com}{\mathbin{{\scriptstyle \circ }}}
\newcommand{\ten}{\mathbin{\otimes}}
\newcommand{\tenll}{\mathbin{\otimes_{R}^{ll}}}

\newcommand{\tenrl}{\mathbin{\otimes_{R}^{rl}}}
\newcommand{\tenrr}{\mathbin{\otimes_{R}^{rr}}}
\newcommand{\tenllM}{\mathbin{\otimes_{\Cc(M)}^{ll}}}

\newcommand{\tenrlM}{\mathbin{\otimes_{\Cc(M)}^{rl}}}

\newcommand{\id}{\mathord{\mathit{id}}}
\newcommand{\pr}{\mathord{\mathit{pr}}}
\newcommand{\inc}{\mathord{\mathit{inc}}}
\newcommand{\supp}{\mathord{\mathrm{supp}}}
\newcommand{\Ker}{\mathord{\mathrm{Ker}}}
\newcommand{\Img}{\mathord{\mathrm{Im}}}
\newcommand{\uni}{\mathord{\mathit{uni}}}
\newcommand{\inv}{\mathord{\mathit{inv}}}
\newcommand{\mlt}{\mathord{\mathit{mlt}}}
\newcommand{\cu}{\mathord{\epsilon}}
\newcommand{\cm}{\mathord{\Delta}}
\newcommand{\Cc}{\mathord{\mathcal{C}^{\infty}_{c}}}
\newcommand{\Cs}[1]{\mathord{\mathcal{C}_{#1}^{\infty}}}
\newcommand{\Esp}{\mathord{\mathcal{E}_{\mathit{sp}}}}
\newcommand{\EspL}{\mathord{\mathcal{E}_{\mathit{sp}}^{\Lambda}}}
\newcommand{\Gsp}{\mathord{\mathcal{G}_{\mathit{sp}}}}
\newcommand{\Sh}{\mathord{\mathrm{Sh}}}
\newcommand{\EtGrpd}{\mathord{\mathrm{EtGrpd}}}
\newcommand{\CoAlg}{\mathord{\mathrm{CoAlg}}}
\newcommand{\HoAlgd}{\mathord{\mathrm{HoAlgd}}}

\title{On spectral representation of coalgebras and Hopf algebroids}
\author{Janez Mr\v{c}un}
\address{Department of Mathematics, University of Ljubljana,
         Jadranska 19, 1000 Ljubljana, Slovenia,
         {\it E-mail:} {\rm janez.mrcun@fmf.uni-lj.si}}
\thanks{This work was supported in part by M\v{S}Z\v{S} grant J1-3148.}

\begin{document}

\begin{abstract}
In this paper we establish a duality between \'{e}tale Lie groupoids
and a class of non-necessarily commutative algebras with a
Hopf algebroid structure. For any \'{e}tale Lie groupoid $G$
over a manifold $M$, the groupoid algebra $\Cc(G)$
of smooth functions with compact support on $G$ has a natural coalgebra structure
over $\Cc(M)$ which makes it into a Hopf algebroid.
Conversely, for any Hopf algebroid $A$ over $\Cc(M)$ we construct the associated
spectral \'{e}tale Lie groupoid $\Gsp(A)$ over $M$ such that $\Gsp(\Cc(G))$
is naturally isomorphic to $G$. Both these constructions are functorial, and
$\Cc$ is fully faithful left adjoint to $\Gsp$.
We give explicit conditions under which a Hopf algebroid is
isomorphic to the Hopf algebroid $\Cc(G)$ of an \'{e}tale Lie groupoid $G$.
We also demonstrate that an analogous duality exists between sheaves on $M$
and a class of coalgebras over $\Cc(M)$.
\end{abstract}

\maketitle

\section*{Introduction}

\'{E}tale Lie groupoids have been successfully used in foliation theory
as suitable geometric models for the transversal structure of foliations.
Many of the important invariants of a foliation reflect
in the structure of the associated holonomy groupoid, which turns into an \'{e}tale
Lie groupoid when restricted to a complete transversal
(see e.g. \cite{Est,Hae58,Hae01}). \'{E}tale Lie groupoids are also used
to represent other geometric objects such as orbifolds and actions of discrete groups.
All these geometric objects can then be studied from the point of view of
non-commutative geometry, by considering the groupoid algebra of
the associated \'{e}tale Lie groupoid, which is the convolution algebra $\Cc(G)$ of
smooth functions with compact support on the \'{e}tale Lie groupoid $G$
introduced by Connes \cite{Con78,Con82,Con94} (see also \cite{CdSW,Ren}).

The aim of this paper is to characterize the algebras which are isomorphic
to groupoid algebras of \'{e}tale Lie groupoids, and to reconstruct the
\'{e}tale Lie groupoid out of its groupoid algebra. We achieve this
in terms of the coalgebra structure on the groupoid algebra
induced by the diagonal map of the groupoid.
If $G$ is an \'{e}tale Lie groupoid over a manifold $M$, then
the coalgebra structure on  $\Cc(G)$ is over the commutative algebra $\Cc(M)$,
and not over the base field. The counit is induced by the target map,
and there is the antipode on $\Cc(G)$ given by the inverse map of $G$.
Such a structure on $\Cc(G)$ is known in the literature as
Hopf algebroid (of quantum groupoid) over $\Cc(M)$ (see \cite{Lu,Mal,Tak,Xu}),
and generalizes the well known Hopf algebra structure on the group algebra associated
to a discrete group. The existence of the Hopf algebroid structure on
$\Cc(G)$ was established in \cite{Mrc99}, but only in the case where
the \'{e}tale Lie groupoid $G$ is Hausdorff. In this paper we extend this definition
to non-Hausdorff groupoids, for which it is essential to take the appropriate
definition of a smooth function with compact support on a non-Hausdorff manifold,
as given in \cite{CraMoe}. It is important to consider the non-Hausdorff groupoids
as well, for instance because the holonomy groupoid of a foliation is often non-Hausdorff.

In the first two sections we concentrate on the coalgebra structure of $\Cc(G)$
over $\Cc(M)$. We show that for any sheaf $E$ on a manifold $M$, the space
$\Cc(E)$ has a natural coalgebra structure over $\Cc(M)$ which is functorial in
$E$. In Section \ref{section:ssc} we then study the structure of an arbitrary coalgebra
$C$ over $\Cc(M)$, and show that for any $x\in M$ there exists a natural quotient
coalgebra
$$ C_{x}=C/N_{x}(C) $$
over the algebra of germs $(\Cs{M})_{x}$ of smooth functions at $x$.
The kernel $N_{x}(C)$ of $C\ra C_{x}$ consists of elements $c\in C$ for which there is
a function $f\in\Cc(M)$ which equals $1$ on a neighbourhood of $x$ and satisfies $fc=0$.
The grouplike elements $G(C_{x})$ of $C_{x}$ are precisely the equivalence classes
of weakly grouplike elements of $C$ normalized at $x$,
i.e. of those elements $c\in C$ for which $\cm(c)=c\ten c'$ for some $c'\in C$
and $\cu(c)$ equals $1$ on a neighbourhood of $x$. We then define the spectral sheaf
$$ \Esp(C) $$
associated to $C$ to be the sheaf on $M$ with the stalk $\Esp(C)_{x}=G(C_{x})$.
We show that this association is also functorial, and that the functor $\Cc$
from the category of sheaves on $M$ to the category of coalgebras over $\Cc(M)$
is fully faithful left adjoint to $\Esp$. In particular, for any sheaf $E$ on $M$
we reconstruct $E$ out of $\Cc(E)$ as the spectral sheaf $\Esp(\Cc(E))$.
We prove that a coalgebra $C$ over $\Cc(M)$ is isomorphic to a sheaf coalgebra $\Cc(E)$
of a sheaf $E$ on $M$ if and only if $C_{x}$ is a free module generated by $G(C_{x})$,
for any $x\in M$. We also give a characterization of the coalgebras
isomorphic to sheaf coalgebras in terms of their weakly grouplike elements.

In Section \ref{section:gha} we show that for any (non-necessarily Hausdorff) \'{e}tale
Lie groupoid $G$ over $M$, the space $\Cc(G)$ of smooth functions with compact support
on $G$ is an algebra with respect to the convolution
product, and that $\Cc(M)$ is a commutative
subalgebra of $\Cc(G)$. Since the target map
$t:G\ra M$  of $G$ is a sheaf on $M$, there is also a natural coalgebra structure
on $\Cc(G)$ over $\Cc(M)$. The inverse map of $G$ induces an antipode $S$ on $\Cc(G)$.
With this structure, $\Cc(G)$ is a Hopf algebroid over $\Cc(M)$.

An arbitrary Hopf algebroid $A$ over $\Cc(M)$ is in particular a coalgebra over $\Cc(M)$,
so we have the spectral sheaf $\Esp(A)$ on $M$ associated to $A$. In Section
\ref{section:sgha} we construct a subsheaf $\Gsp(A)$ of $\Esp(A)$.
The stalk $\Gsp(A)_{y}$ of
$\Gsp(A)$ at any $y\in M$ is the subset of $G(A_{y})$ of those grouplike elements of
$A_{y}$ which are equivalence classes of $S$-invariant weakly grouplike elements of $A$,
i.e. elements $a$ of $A$ for which $\cm(a)=a\ten a'$ and $\cm(S(a))=S(a')\ten S(a)$
for some $a'\in A$. We then show that $\Gsp(A)$ has a natural structure of an \'{e}tale
Lie groupoid over $M$, with the sheaf projection $\Gsp(A)\ra M$ as the target map.
The constructed \'{e}tale Lie groupoid
$$ \Gsp(A) $$
over $M$ is referred to as spectral \'{e}tale groupoid associated to $A$.
The functor $\Cc$ from the category of \'{e}tale Lie groupoids over $M$
to the category of Hopf algebroids over $\Cc(M)$
is fully faithful left adjoint to $\Gsp$, and $\Gsp(\Cc(G))$ is naturally
isomorphic to $G$, for any \'{e}tale Lie groupoid $G$.
We show that a Hopf algebroid $A$ over $\Cc(M)$ is isomorphic to a groupoid
Hopf algebroid $\Cc(G)$ of an \'{e}tale Lie groupoid $G$ over $M$ if and only if
$A_{y}$ is a free module generated by $\Gsp(A)_{y}$,
for any $y\in M$. We also express this condition
in terms of the $S$-invariant weakly grouplike elements of $A$.

Throughout this paper, all manifolds and maps between them are assumed to be smooth,
i.e. of the class $\mathcal{C}^{\infty}$.
However, all the results and arguments hold true if one replaces this by any class
of differentiability $\mathcal{C}^{k}$, $k=0,1,2,\ldots\,$.
In this paper, we do not consider the generalized maps between sheaves and \'{e}tale
Lie groupoids over different bases, and associated Morita theory. In some extend we
have done this in \cite{Mrc01} for Hausdorff \'{e}tale Lie groupoids only, however
non-Hausdorff \'{e}tale Lie groupoids can be treated along the same lines as well.

\section{Sheaf coalgebras}  \label{section:sc}

The aim of this section is to show that the space of smooth
functions with compact support, defined on the total space
of a sheaf on a Hausdorff manifold $M$,
has a natural coalgebra structure over the ring $\Cc(M)$ of smooth
functions with compact support on $M$. First, let
us briefly recall the notion of a coalgebra over a ring,
to fix the notations (for details, see e.g. \cite{MilMoo}).

Note that the ring $\Cc(M)$ has the identity
only if $M$ is compact, but it always has local identities.
Recall that a commutative ring $R$ has {\em local identities} if for
any $r_{1},\ldots,r_{k}\in R$
there exists $r\in R$ such that $rr_{i}=r_{i}$ for any $i=1,\ldots,k$.
More generally, we say that a (left) $R$-module $C$ is {\em locally unitary}
if for any $c_{1},\ldots,c_{k}\in C$ there exists $r\in R$ such that
$rc_{i}=c_{i}$ for any $i=1,\ldots,k$.

Let $R$ be a ring with local identities and let $C$ be a
locally unitary (left) $R$-module.
A \textit{coalgebra structure} on $C$ over $R$ consists of
two $R$-linear maps, the {\em comultiplication}
$$ \cm:C\lra C\ten_{R} C $$
and the {\em counit}
$$ \cu:C\lra R \;,$$
which satisfy $(\id\ten\cu)\com\cm=\id$
and $(\cu\ten\id)\com\cm=\id$. Here we identified the modules
$C$, $R\ten_{R}C$ and $C\ten_{R}R$ in the natural way,
which is possible because $C$ is locally unitary.
An $R$-coalgebra (or a coalgebra over $R$)
is a locally unitary $R$-module $C$ equipped with
a coalgebra structure $(\cm,\cu)$ over $R$. We say that
an $R$-coalgebra $C=(C,\cm,\cu)$ is \textit{coassociative}
if $(\cm\ten\id)\com\cm=(\id\ten\cm)\com\cm$. We say that $C$ is
\textit{cocommutative} if $\sigma\com\cm=\cm$,
where $\sigma:C\ten_{R} C\ra C\ten_{R} C$ is the flip between the factors, i.e.
$\sigma(c\ten c')=c'\ten c$.
If $D$ is another $R$-coalgebra, then a map $\alpha:C\ra D$
is a homomorphism of $R$-coalgebras if $\alpha$ is $R$-linear,
$(\alpha\ten\alpha)\com\cm=\cm\com\alpha$ and $\cu=\cu\com\alpha$.
Denote by
$$ \CoAlg(R) $$
the category of $R$-coalgebras and homomorphisms between them.

Let $R$ be a commutative ring with local identities and let
$C$ be a coalgebra over $R$.
We say that an element $c\in C$ is {\em weakly grouplike} if
$\cm(c)=c\ten c'$ for some $c'\in C$. Denote by
$$ G_{w}(C) $$
the set of all weakly grouplike elements of $C$.
If $R$ has an identity, the set of {\em grouplike elements} of $C$ is defined
(c.f. \cite{MilMoo}) as
$$ G(C)=\{c\in C\,|\,\cm(c)=c\ten c\,,\;\cu(c)=1\} \;.$$
Note that for a coalgebra $C$ over a field $\FF$,
the weakly grouplike elements of $C$ are exactly all the products $rc$ of
any grouplike element $c\in C$ with any $r\in\FF$.
\vspace{2mm}

\Rem
Note that if $C$ is generated by $G_{w}(C)$ as an $R$-module, then
$C$ is cocommutative and coassociative. Indeed, if $c\in G_{w}(C)$ with
$\cm(c)=c\ten c'$, then $c=\cu(c)c'=\cu(c')c$,
$\cm(c)=\cu(c)c'\ten c'=c'\ten c$ and
$(\id\ten\cm)(\cm(c))=\cu(c)c'\ten c'\ten c'=(\cm\ten\id)(\cm(c))$.
Since $G_{w}(C)$ is clearly closed under the $R$-action, it follows that
$C$ is generated by $G_{w}(C)$ as an $R$-module if and only if it
is generated by $G_{w}(C)$ as an abelian group.
\vspace{2mm}

It is well known than grouplike elements in a coalgebra over a field
are automatically linearly independent, but this is unfortunately no longer true
if we are working over a commutative ring with divisors of zero (and identity).
We give here sufficient conditions under which this is true,
and an example which illustrates that the conditions are necessary. We will use
this result in the following sections.

\begin{prop}  \label{prop1}
Let $R$ be a commutative ring with identity without nilpotents of order $2$.
Let $C$ be a coalgebra over $R$, and let $L$ be a subset of $G(C)$.
If $C$ is a flat $R$-module and if any two grouplike elements in $L$
are $R$-linearly independent,
then $L$ is $R$-linearly independent.
\end{prop}
\begin{proof}
We have to show that if
$$ r_{1}c_{1}+\ldots+r_{k}c_{k}=0 $$
for some grouplike elements $c_{1},\ldots,c_{k}\in L$ and
$r_{1},\ldots,r_{k}\in R$, then $r_{i}=0$ for any $1\leq i\leq k$.
We prove this by induction on $k$. If $k=1$ we get $r_{1}=0$
because $\cu(c_{1})=1$, while for $k=2$
the statement holds by assumption. So let $k\geq 3$, and take any
$1\leq l\leq k$. By the induction hypothesis we know
that $c_{1},\ldots,c_{l-1},c_{l+1},\ldots,c_{k}$ are $R$-linearly independent.
The flatness of $C$ ensures that
$$ \{c_{i}\ten c_{j}\,|\, 1\leq i,j\leq k\,,\;i\neq l\,,\;j\neq l \} $$
is $R$-linearly independent in $C\ten_{R}C$.
Now
\begin{eqnarray*}
0
& = & r_{l}\cm(\sum_{i}r_{i}c_{i})= r_{l}\sum_{i}r_{i}c_{i}\ten c_{i} \\
& = & \sum_{i\neq l}r_{l}r_{i}c_{i}\ten c_{i}+r_{l}c_{l}\ten r_{l}c_{l} \\
& = & \sum_{i\neq l}r_{l}r_{i}c_{i}\ten c_{i}
      +(-\sum_{i\neq l}r_{i}c_{i})\ten (-\sum_{j\neq l}r_{j}c_{j}) \\
& = & \sum_{i\neq l}r_{l}r_{i}c_{i}\ten c_{i}
      +\!\!\!\sum_{i\neq l, j\neq l}\!\!\!r_{i}r_{j}c_{i}\ten c_{j}\;.
\end{eqnarray*}
Thus for any $i\neq l$ we have $r_{l}r_{i}+r_{i}r_{i}=0$, and furthermore
$r_{i}r_{j}=0$ for any $j\neq l$, $j\neq i$. Since this is true for any
$1\leq l\leq k$, we can conclude that for any $1\leq i\leq k$ we have
$r_{i}r_{i}=0$ and thus $r_{i}=0$.
\end{proof}

\begin{ex}  \rm \label{ex2}
Let $R$ be a commutative ring with identity and with divisors of zero.
Choose non-zero $p,q\in R$ with $pq=0$. Let $C$ be the
free $R$-module
on two generators $e_{1}$ and $e_{2}$, and define an $R$-coalgebra structure
on $C$ by
$$ \cm(e_{1})=e_{1}\ten e_{1}\;,\;\;\;\;\;\;\;\;\;\; \cu(e_{1})=1\;, $$
$$ \cm(e_{2})=e_{1}\ten e_{2}+e_{2}\ten e_{1}+p e_{2}\ten e_{2}\,,
   \;\;\;\;\;\;\;\;\;\;\cu(e_{2})=0 \;.$$
One can easily check that this is indeed a coassociative cocommutative
coalgebra structure on $C$. One grouplike element of $C$ is clearly $e_{1}$,
while another is
$$ c=e_{1}+pe_{2} \;.$$
But these two are not $R$-linearly independent because $qe_{1}=qc$.
\end{ex}

We shall now describe the coalgebra of smooth functions on a sheaf on
a Hausdorff manifold.
For this to be correct, one should be careful with the notion
of a smooth function with compact support on a non-Hausdorff manifold.
It turns out that the definition given in \cite{CraMoe},
which we recall bellow, is the most suitable for our purpose.

Let $E$ be a non-necessarily Hausdorff manifold.
Denote by $\Cs{E}$ the sheaf of germs of smooth
functions on $E$ with values in $\FF\in\{\RR,\CC\}$.
The stalk $(\Cs{E})_{y}$ of this
sheaf at any point $y\in E$ is a commutative algebra with identity.
If $U$ is a Hausdorff open subset of $E$ then the continuous sections
of $\Cs{E}$ with compact support in $U$ are just the
smooth functions with compact support in $U$,
$$ \Gamma_{c}(U,\Cs{E})= \Cc(U)\;. $$
The extension of such a section by zero to all of $E$
may not be continuous if $E$ is not Hausdorff.
Therefore, extension by zero gives a monomorphism
$$ \Cc(U)\lra\Gamma_{\delta}(E,\Cs{E})\;,$$
where $\Gamma_{\delta}(E,\Cs{E})$ denotes the space of all
non-necessarily continuous global sections of the sheaf $\Cs{E}$.
For any $f\in\Cc(U)$ we shall denote the corresponding extension by zero
of $f$ again by $f\in \Gamma_{\delta}(E,\Cs{E})$, or by
$f|_{U}^{E}$ if necessary.
For any section $u\in \Gamma_{\delta}(E,\Cs{E})$
and any $y\in E$, denote by
$$ u_{y}\in(\Cs{E})_{y} $$
the value of this section at $y$, and by
$$ u(y)=u_{y}(y)\in\FF $$
the value of the germ $u_{y}$ at $y$.

The space $\Cc(E)$ of {\em smooth functions with compact support} on $E$
is the image of the map $\bigoplus_{U}\Cc(U)\ra\Gamma_{\delta}(E,\Cs{E})$,
where $U$ ranges over all Hausdorff open subsets of $E$,
$$ \Cc(E)=\Img(\bigoplus_{U}\Cc(U)\lra\Gamma_{\delta}(E,\Cs{E}))\;.$$
One can easily see, by using partitions of unity, that the map
$$  \bigoplus_{i\in I}\Cc(U_{i})\lra\Cc(E) $$
is surjective for any Hausdorff open cover $(U_{i})_{i\in I}$ of $E$.
In particular, this definition agrees with the classical one if $E$ is
Hausdorff. For any $u\in\Cc(E)$ define the \textit{support of} $u$ by
$$ \supp(u)=\{y\in E\,|\,u_{y}\neq 0\}\;.$$
Observe that
if $\supp(u)\subset U$ for some Hausdorff open subset $U$ of $E$, then
$u=f|_{U}^{E}$ for a unique $f\in\Cc(U)$. In this case, the support of
$u$ equals the support of $f$ in $U$ defined in the classical sense.
Note, however, that $\Cc(E)$ is closed for the pointwise multiplication
only if $E$ is Hausdorff.

Let $F$ be another non-necessarily Hausdorff manifold, and
let $\phi:E\ra F$ be a smooth map. Then for any $y\in E$ the composition
with $\phi$ induces a homomorphism of commutative algebras
$$ \phi_{y}^{\ast}:(\Cs{F})_{\phi(y)}\lra (\Cs{E})_{y} \;.$$

Suppose that $\phi:E\ra F$ is a local diffeomorphism.
Then $\phi_{y}^{\ast}$ is an isomorphism, and we denote its inverse by
$$ \phi_{\ast y}=(\phi_{y}^{\ast})^{-1}\;.$$
An open subset $U$ of $E$ will be called $\phi$-\textit{elementary} if
$\phi|_{U}$ is injective. Note that Hausdorff $\phi$-elementary
subsets of $E$ form a basis of $E$.
We can define a linear map
$$ \phi_{+}:\Cc(E)\lra\Cc(F) $$
by
$$ \phi_{+}(u)_{z}=\!\!\!\sum_{y\in \phi^{-1}(z)}
   \!\!\!\phi_{\ast y}(u_{y})\;. $$
One can easily see that $\phi_{+}(u)$ is a smooth function with
compact support on $F$. Indeed, one can assume without loss of generality
that the support of $u$ lies in
a $\phi$-elementary Hausdorff subset $U$ of $E$. Then $V=\phi(U)$ is
an open Hausdorff subset of $F$, $\phi|_{U}:U\ra V$ is a diffeomorphism
and $\phi_{+}(u)=u\com(\phi|_{U})^{-1}$. In this way,
we extended $\Cc$ to a functor defined on local diffeomorphisms
between non-necessarily Hausdorff smooth manifolds.

Let $M$ be a Hausdorff smooth manifold.
Recall that a \textit{sheaf} (of sets) on $M$ is
a local diffeomorphism
$$ \pi:E\lra M\;,$$
where $E$ is a non-necessarily Hausdorff manifold, called the
{\em total space} of the sheaf. The sheaf $\pi:E\ra M$ is often denoted
simply by $E$, while $\pi$ is referred to as the
{\em projection} of the sheaf $E$.
The {\em stalk} of $E$ at $x\in M$ is the fiber $E_{x}=\pi^{-1}(x)$.
A morphism from a sheaf $\pi:E\ra M$ to another sheaf
$\xi:F\ra M$ over $M$ is a smooth map $\phi:E\ra F$ over $M$, i.e. a map
which satisfies $\pi=\xi\com\phi$.
Any such map is necessarily a local diffeomorphism.
The category of sheaves on $M$ and morphisms between them
is denoted by
$$\Sh(M)\;.$$

Let $M$ be a Hausdorff manifold and $\pi:E\ra M$ a sheaf on $M$.
First note that any $\pi$-elementary subset $U$ of
$E$ is Hausdorff.
If $U\subset E$ is $\pi$-elementary and $f\in\Cc(\pi(U))$,
we have
$$ f\com\pi|_{U}=(f\com\pi|_{U})|_{U}^{E}\in\Cc(E)\;.$$
Any element of $\Cc(E)$ can be written as a sum of elements
of this form.
There is a natural (left) $\Cc(M)$-module structure on
$\Gamma_{\delta}(E,\Cs{E})$, given
for any $p\in\Cc(M)$ and $u\in\Gamma_{\delta}(E,\Cs{E})$ by
$$ (pu)_{y}=\pi^{\ast}_{y}(p_{\pi(y)})u_{y} \;,$$
which restricts to a $\Cc(M)$-module structure on $\Cc(E)$.
In fact, we have
$$ p(f\com\pi|_{U})=pf\com\pi|_{U}\;.$$
If $\xi:F\ra M$ is another sheaf on $M$ and $\phi:E\ra F$ a
morphism of sheaves over $M$, then $\phi_{+}:\Cc(E)\ra\Cc(F)$
is a homomorphism of $\Cc(M)$-modules.
If $U\subset E$ is $\pi$-elementary, then
$\phi(U)$ is a $\xi$-elementary subset of $F$ and
$$ \phi_{+}(f\com\pi|_{U})=f\com\xi|_{\phi(U)} $$
for any $f\in\Cc(\pi(U))=\Cc(\xi(\phi(U)))$.
Therefore $\Cc$ can be regarded as
a functor from the category of sheaves on $M$ to the category
of $R$-modules over $\Cc(M)$.

To define a coalgebra structure on $\Cc(E)$ we essentially need the
following fact:

\begin{prop}  \label{prop3}
Let $M$ be a Hausdorff manifold, and
let $\pi:E\ra M$ and $\xi:F\ra M$ be two sheaves on $M$.
Denote by $E\times_{M}F$ the
fibered product of $E$ and $F$, as in the following diagram:
$$
\xymatrix{
E\times_{M}F \ar[d]_{\pr_{1}} \ar[r]^-{\pr_{2}} & F \ar[d]^{\xi} \\
E \ar[r]^{\pi} & M
}
$$
Then $\pi\com\pr_{1}=\xi\com\pr_{2}:E\times_{M}F\ra M$ is a sheaf
on $M$, and
there is a natural isomorphism of $\Cc(M)$-modules
$$ \Omega_{\pi,\xi}:\Cc(E)\ten_{\Cc(M)}\Cc(F)\lra \Cc(E\times_{M}F) $$
given by
$$ \Omega_{\pi,\xi}(u\ten v)_{(y,z)}
   =(\pr_{1})^{\ast}_{(y,z)}(u_{y}) (\pr_{2})^{\ast}_{(y,z)}(v_{z}) \;,$$
for any $u\in\Cc(E)$, $v\in\Cc(F)$ and $(y,z)\in E\times_{M}F$.
If $U\subset E$ is $\pi$-elementary and $V\subset F$ is
$\xi$-elementary, then $U\times_{M}V$ is a $(\pi\com\pr_{1})$-elementary
subset of $E\times_{M}F$ and
$$ \Omega_{\pi,\xi}(f\com\pi|_{U}\ten p\com\xi|_{V})
   =fp\com (\pi\com\pr_{1})|_{U\times_{M}V} $$
for any $f\in\Cc(\pi(U))$ and $p\in\Cc(\xi(V))$.
\end{prop}
\begin{proof}
First, the maps $\pr_{1}$ and $\pr_{2}$ are local diffeomorphisms
because they are pull backs of local diffeomorphisms.
Therefore $\omega=\pi\com\pr_{1}=\xi\com\pr_{2}$ is a local diffeomorphism
as well. It is clear that
$\Omega_{\pi,\xi}(u\ten v)$ is section in
$\Gamma_{\delta}(E\times_{M}F,\Cs{E\times_{M}F})$. The equation
$$ \Omega_{\pi,\xi}(f\com\pi|_{U}\ten p\com\xi|_{V})
   =fp\com\omega|_{U\times_{M}V} $$
follows directly from the definition, and implies that
$\Omega_{\pi,\xi}(u\ten v)$ is indeed a smooth function with compact
support on $E\times_{M}F$. One can also easily check that
$\Omega_{\pi,\xi}$ is well defined on the tensor product over $\Cc(M)$ and
$\Cc(M)$-linear.

To show that $\Omega_{\pi,\xi}$ is surjective it is sufficient to prove that
$f\com\omega|_{W}$ is in the image of $\Omega_{\pi,\xi}$ for any
$\omega$-elementary subset
$W$ of $E\times_{M}F$ and any $f\in\Cc(\omega(W))$.
To this end, note that $\pr_{1}(W)$ is $\pi$-elementary,
$\pr_{2}(W)$ is $\xi$-elementary and $W=\pr_{1}(W)\times_{M}\pr_{2}(W)$.
Choose any function $p\in\Cc(\omega(W))$ with $fp=f$. Now we have
$$ f\com\omega|_{W}
   =\Omega_{\pi,\xi}(f\com\pi|_{\pr_{1}(W)}\ten p\com\xi|_{\pr_{2}(W)})\;.$$

Next, we have to show that $\Omega_{\pi,\xi}$ is injective.
Take any element in the kernel of $\Omega_{\pi,\xi}$
and write it in the form
$$ \sum_{i=1}^{k} f_{i}\com\pi|_{U_{i}} \ten p_{i}\com\xi|_{V_{i}}\;,$$
for some $\pi$-elementary subsets $U_{i}$ of $E$,
$\xi$-elementary subsets $V_{i}$ of $F$,
$f_{i}\in\Cc(\pi(U_{i}))$ and $p_{i}\in\Cc(\xi(V_{i}))$, $i=1,\ldots,k$.
We assume that $k$ is the smallest number for which this can be done.
We shall now show that
the assumption $k>0$ leads to a contradiction, and this will conclude our proof.

First note that we can assume without loss of generality that
$\pi(U_{i})=\xi(V_{i})$ and that
$f_{i}p_{i}=f_{i}$,
for any $i=1,\ldots,k$. Indeed,
we can choose functions $\eta_{i}\in\Cc(\pi(U_{i})\cap\xi(V_{i}))$ and
$\theta_{i}\in\Cc(\xi(V_{i}))$ such that
$\eta_{i}(f_{i}p_{i})=f_{i}p_{i}$ and
$p_{i}\theta_{i}=p_{i}$.
Put $U'_{i}=(\pi|_{U_{i}})^{-1}(\xi(V_{i}))$ and
$V'_{i}=(\xi|_{V_{i}})^{-1}(\pi(U_{i}))$.
Then $p_{i}f_{i}\com\pi|_{U'_{i}} \ten \eta_{i}\theta_{i}\com\xi|_{V'_{i}}$
satisfy the required condition and equals
\begin{eqnarray*}
p_{i}f_{i}\com\pi|_{U'_{i}} \ten \eta_{i}\theta_{i}\com\xi|_{V'_{i}}
   & = & \eta_{i}p_{i}f_{i}\com\pi|_{U'_{i}} \ten \theta_{i}\com\xi|_{V_{i}} \\
   & = & p_{i}f_{i}\com\pi|_{U'_{i}} \ten \theta_{i}\com\xi|_{V_{i}} \\
   & = & f_{i}\com\pi|_{U_{i}} \ten p_{i}\theta_{i}\com\xi|_{V_{i}} \\
   & = & f_{i}\com\pi|_{U_{i}} \ten p_{i}\com\xi|_{V_{i}}\;.
\end{eqnarray*}

Put $W_{i}=U_{i}\times_{M}V_{i}$,  for any $i=1,\ldots,k$.
By our assumption we have $\pi(U_{i})=\xi(V_{i})=\omega(W_{i})$.
Since $\sum_{i=1}^{k} f_{i}\com\pi|_{U_{i}} \ten p_{i}\com\xi|_{V_{i}}$
is in the kernel of $\Omega_{\pi,\xi}$, we have
$$ \sum_{i=1}^{k} p_{i}f_{i}\com\omega|_{W_{i}}
   =\sum_{i=1}^{k} f_{i}\com\omega|_{W_{i}}=0\;.$$
Now let $X_{i}=W_{i}\cap W_{k}$, for $i=1,\ldots,k-1$. We claim that
$$ \supp(f_{k}\com\omega|_{W_{k}})
   \subset \bigcup_{i=1}^{k-1}X_{i}\subset W_{k}\;.$$
Indeed, if $(y,z)\in\supp(f_{k}\com\omega|_{W_{k}})$ then
there exists
$1\leq i\leq k-1$ with $(y,z)\in\supp(f_{i}\com\omega|_{W_{i}})\subset W_{i}$
because $\sum_{i=1}^{k} f_{i}\com\omega|_{W_{i}}=0$. Therefore
$f_{k}\in\Cc(\bigcup_{i=1}^{k-1}\omega(X_{i}))$, hence we can choose
$\rho\in\Cc(\bigcup_{i=1}^{k-1}\omega(X_{i}))$
with $\rho f_{k}=f_{k}$. Furthermore, we can write
$\rho=\sum_{i=1}^{k-1}\rho_{i}$ for some
$\rho_{i}\in\Cc(\omega(X_{i}))$, and we can choose
$\nu_{i}\in\Cc(\omega(W_{i}))$ such that
$\nu_{i}f_{i}=f_{i}$ and  $\nu_{i}\rho_{i}=\rho_{i}$, for $i=1,\ldots,k-1$.
Now we first observe that
\begin{eqnarray*}
f_{i}\com\pi|_{U_{i}} \ten p_{i}\com\xi|_{V_{i}}
   & = & \nu_{i}f_{i}\com\pi|_{U_{i}} \ten p_{i}\com\xi|_{V_{i}} \\
   & = & f_{i}\com\pi|_{U_{i}} \ten \nu_{i}p_{i}\com\xi|_{V_{i}} \\
   & = & p_{i}f_{i}\com\pi|_{U_{i}} \ten \nu_{i}\com\xi|_{V_{i}} \\
   & = & f_{i}\com\pi|_{U_{i}} \ten \nu_{i}\com\xi|_{V_{i}}
\end{eqnarray*}
and
\begin{eqnarray*}
\rho_{i}f_{k}\com\pi|_{U_{k}} \ten p_{k}\com\xi|_{V_{k}}
& = & f_{k}\com\pi|_{U_{k}} \ten \nu_{i}\rho_{i}p_{k}\com\xi|_{V_{k}} \\
& = & f_{k}\com\pi|_{U_{k}} \ten \nu_{i}\rho_{i}p_{k}\com\xi|_{\pr_{2}(X_{i})} \\
& = & \rho_{i}p_{k}f_{k}\com\pi|_{U_{k}} \ten \nu_{i}\com\xi|_{V_{i}} \\
& = & \rho_{i}p_{k}f_{k}\com\pi|_{\pr_{1}(X_{i})} \ten \nu_{i}\com\xi|_{V_{i}} \\
& = & \rho_{i}f_{k}\com\pi|_{U_{i}} \ten \nu_{i}\com\xi|_{V_{i}}\;.
\end{eqnarray*}
Therefore we have
\begin{eqnarray*}
\lefteqn{\sum_{i=1}^{k} f_{i}\com\pi|_{U_{i}} \ten p_{i}\com\xi|_{V_{i}}}
&   &  \\
& = &  \sum_{i=1}^{k-1} f_{i}\com\pi|_{U_{i}} \ten p_{i}\com\xi|_{V_{i}}
  + \sum_{i=1}^{k-1} \rho_{i}f_{k}\com\pi|_{U_{k}} \ten p_{k}\com\xi|_{V_{k}}  \\
& = &  \sum_{i=1}^{k-1} f_{i}\com\pi|_{U_{i}} \ten \nu_{i}\com\xi|_{V_{i}}
  + \sum_{i=1}^{k-1} \rho_{i}f_{k}\com\pi|_{U_{i}}\ten \nu_{i}\com\xi|_{V_{i}} \\
& = &  \sum_{i=1}^{k-1} (f_{i}
  +\rho_{i}f_{k})\com\pi|_{U_{i}} \ten \nu_{i}\com\xi|_{V_{i}}\;.
\end{eqnarray*}
Thus we reduced the number of summands for
one, which contradicts our assumption.

Finally, we will check that $\Omega_{\pi,\xi}$ is natural in $\pi$ and $\xi$.
Suppose that $\pi':E'\ra M$ and $\xi':F'\ra M$ are another two sheaves on $M$ and
that $\phi:E\ra E'$ and $\psi:F\ra F'$ are morphisms of sheaves over $M$.
If $U\subset E$ is $\pi$-elementary and $V\subset F$ is
$\xi$-elementary, then $\phi(U)\subset E'$ is $\pi'$-elementary
and $\psi(V)\subset F'$ is $\xi'$-elementary.
For any $f\in\Cc(\pi(U))$ and $p\in\Cc(\xi(V))$ we have
\begin{eqnarray*}
(\phi\times\psi)_{+}(\Omega_{\pi,\xi}(f\com\pi|_{U}\ten p\com\xi|_{V}))
& = & (\phi\times\psi)_{+}(fp\com (\pi\com\pr_{1})|_{U\times_{M}V}) \\
& = & fp\com(\pi'\com\pr_{1})|_{\phi(U)\times_{M}\psi(V)} \\
& = & \Omega_{\pi',\xi'}(f\com\pi'|_{\phi(U)}\ten p\com\xi'|_{\psi(V)}) \\
& = & \Omega_{\pi',\xi'}((\phi_{+}\ten\psi_{+})(f\com\pi|_{U}\ten p\com\xi|_{V})) \;,
\end{eqnarray*}
hence
$(\phi\times\psi)_{+}\com\Omega_{\pi,\xi}=\Omega_{\pi',\xi'}\com (\phi_{+}\ten\psi_{+})$.
\end{proof}

Let $\pi:E\ra M$ be a sheaf on $M$.
We define a coalgebra structure on $\Cc(E)$ over
$\Cc(M)$ as follows:
the comultiplication $\cm$ is given by the composition
$$
\xymatrix{
\Cc(E) \ar[r]^-{\mathrm{diag}_{+}} &  \Cc(E\times_{M}E) \ar[r]^-{\Omega_{\pi,\pi}^{-1}}
  & \Cc(E)\ten_{\Cc(M)} \Cc(E) \;,
}
$$
where $\mathrm{diag}:E\ra E\times_{M}E$ is the diagonal map defined with
$\mathrm{diag}(y)=(y,y)$, while the counit is
$$ \cu=\pi_{+}:\Cc(E)\lra\Cc(M)\;.$$
If $U$ is any $\pi$-elementary subset of $E$ and $f\in\Cc(\pi(U))$, we have
$$ \cm(f\com \pi|_{U})=(f\com \pi|_{U})\ten (\eta\com \pi|_{U})
   =(\eta\com \pi|_{U})\ten (f\com \pi|_{U}) \;,$$
where $\eta$ is any element of $\Cc(\pi(U))$ with $f\eta=f$. Furthermore,
$$ \cu(f\com \pi|_{U})=f\;.$$
One can now easily check that $\Cc(E)=(\Cc(E),\cm,\cu)$ is a cocommutative
coassociative coalgebra over $\Cc(M)$.

Suppose that $\xi:F\ra M$ is another sheaf on $M$, and
let $\phi:E\ra F$ be morphism of sheaves on $M$.
Then $\phi_{+}:\Cc(E)\ra\Cc(F)$ is a homomorphism of $\Cc(M)$-coalgebras.
Indeed, we have
\begin{eqnarray*}
\cm\com\phi_{+}
& = & \Omega_{\xi,\xi}^{-1}\com\mathrm{diag}_{+}\com\phi_{+}
      =\Omega_{\xi,\xi}^{-1}\com(\phi\times\phi)_{+}\com\mathrm{diag}_{+} \\
& = & (\phi_{+}\ten\phi_{+})\com\Omega_{\pi,\pi}^{-1}\com\mathrm{diag}_{+}
      =(\phi_{+}\ten\phi_{+})\com\cm\;,
\end{eqnarray*}
and
$$ \cu\com\phi_{+}=\xi_{+}\com\phi_{+}=(\xi\com\phi)_{+}=\pi_{+}=\cu\;.$$
Therefore we may conclude that we have a functor
$$ \Cc:\Sh(M)\lra\CoAlg(\Cc(M))\;.$$
The next section will show that this functor is fully faithful
and has a right adjoint.

We can explicitly describe weakly grouplike elements in a sheaf coalgebra as follows:

\begin{prop} \label{prop4}
Let $M$ be a Hausdorff manifold and let $\pi:E\ra M$ be a sheaf on $M$.
An element $u\in\Cc(E)$ is weakly grouplike if and only if
$u\in\Cc(U)$ for some  $\pi$-elementary subset $U$ of $E$.
\end{prop}
\begin{proof}
Let $u\in\Cc(E)$ be weakly grouplike, with
$\cm(u)=u\ten u'$  for some $u'\in\Cc(E)$.
This equality implies that for any
$(y,z)\in E\times_{M}E$, $y\neq z$, we have
$$ (\pr_{1})^{\ast}_{(y,z)}(u_{y}) (\pr_{2})^{\ast}_{(y,z)}(u'_{z})=0 $$
and
$$ (\pr_{1})^{\ast}_{(y,y)}(u_{y})
   (\pr_{2})^{\ast}_{(y,y)}(u'_{y})=(\pr_{1})^{\ast}_{(y,y)}(u_{y}) \;.$$
So if $u_{y}\neq 0$ then $u'(z)=0$ and $u'(y)=1$, hence $u_{z}=0$. In particular,
the map $\pi$ is injective on $\supp(u)$.

This yields that for any $y\in\supp(u)$
we can choose a $\pi$-elementary subset $V$ of $E$
with $y\in V$ such that
$V\cap\supp(u)=\pi^{-1}(\pi(V))\cap \supp(u)$. Indeed,
write $u=\sum_{i=1}^{k}f_{i}\com \pi|_{U_{i}}$ for some
$\pi$-elementary $U_{i}\subset E$ and $f_{i}\in\Cc(\pi(U_{i}))$.
Let $y_{0},y_{1},\ldots,y_{n}$ be all elements of
$\pi^{-1}(\pi(y))\cap\bigcup_{i=1}^{k} U_{i}$, with $y_{0}=y$. For any
$0\leq l\leq n$, let $I_{l}$ denote the set of all $1\leq i\leq k$
with $y_{l}\in U_{i}$. If $l\geq 1$ we have
$u_{y_{l}}=0$ and hence
$$ \sum_{i\in I_{l}}(f_{i})_{\pi(y)}=0\;.$$
Therefore we can choose a neighbourhood $W$ of $\pi(y)$ so small that
$W\subset \pi(\bigcap_{i\in I_{l}}U_{i})$
for any $0\leq l\leq n$, and
$$ \sum_{i\in I_{l}} f_{i}|_{W}=0 $$
for any $l\geq 1$. Furthermore we can choose $W$ so small that
$W\cap \supp(f_{i})=\emptyset$ for any $i$ with $\pi(y)\not\in \pi(U_{i})$.
Now $V=\pi^{-1}(W)\cap \bigcap_{i\in I_{0}}U_{i}$
is $\pi$-elementary satisfying
$V\cap\supp(u)= \pi^{-1}(\pi(V))\cap\supp(u)$.

Since $\supp(\cu(u))=\pi(\supp(u))$ is compact, we can choose
a finite number of $\pi$-elementary subsets $V_{1},\ldots,V_{m}$
of $E$ which all satisfy
$V_{j}\cap\supp(u)=\pi^{-1}(\pi(V_{j}))\cap \supp(u)$, such that
$\pi(V_{1}),\ldots,\pi(V_{m})$ cover $\supp(\cu(u))$. This yields
that $V_{1},\ldots,V_{m}$ cover $\supp(u)$ as well.
Choose $\rho\in\Cc(\bigcup_{j=1}^{m}\pi(V_{j}))$
with $\rho\cu(u)=\cu(u)$, and write $\rho=\sum_{j=1}^{m}\rho_{j}$ for some
$\rho_{j}\in\Cc(\pi(V_{j}))$. This implies that $\rho u=u$ and
that $\rho_{j}u\in\Cc(V_{j})$.
Choose compact sets $K_{j}\subset V_{j}$ and open sets $Z_{j}\subset K_{j}$
such that $\supp(\rho_{j}u)\subset Z_{j}$,
for any $1\leq j\leq m$. Note that $Z=\bigcup_{j=1}^{m}Z_{j}$
is an open neighbourhood of $\supp(u)$ because $\rho u=u$.

Now for any $y\in\supp(u)$ we can choose an open neighbourhood $Y_{y}\subset Z$
such that for any $1\leq j\leq m$
we have either $Y_{y}\subset V_{j}$ if $y\in V_{j}$,
or $\pi(Y_{y})\cap \pi(K_{j})=\emptyset$ if $y\not\in V_{j}$.
Indeed, if $y\not\in V_{j}$ then we also have
$\pi(y)\not\in\pi(V_{j})$ because
$V_{j}\cap\supp(u)=\pi^{-1}(\pi(V_{j}))\cap \supp(u)$. Then
$$ Y=\!\!\!\bigcup_{y\in\supp(u)}\!\!\!Y_{y}\subset Z $$
is an open neighbourhood of $\supp(u)$
and is $\pi$-elementary. To see this, suppose that $z,w\in Y$ with
$\pi(z)=\pi(w)$. Then there exists $y\in\supp(u)$ with $z\in Y_{y}$ and there
exists $1\leq j\leq m$ such that $w\in Z_{j}\subset K_{j}$. If
$y\not\in V_{j}$ then $\pi(z)\in \pi(Y_{y})\cap \pi(K_{j})=\emptyset$,
a contradiction.
Therefore we must have $y\in V_{j}$, so $Y_{y}\subset V_{j}$ and
hence $z,w\in V_{j}$. Since $\pi(z)=\pi(w)$ and $V_{j}$ is $\pi$-elementary,
it follows that $z=w$.
\end{proof}

\section{Spectral sheaf of a coalgebra}  \label{section:ssc}

Let $M$ be a Hausdorff manifold.
Throughout this section we shall write $R=\Cc(M)$.
For any left $R$-module $C$ and any $x\in M$ we denote
$$ N_{x}(C)=\{c\in C\,|\, fc=0 \;\textrm{for some}\; f\in\Cc(M)\;
   \mathrm{with}\; f_{x}=1\}\;. $$
Note that $N_{x}(C)$ is an $R$-submodule of $C$.
The quotient of $C$ over $N_{x}(C)$ is also an $R$-module, and we shall
denote this quotient by
$$ C_{x}=C/N_{x}(C)\;.$$
For any $c\in C$ we shall  write
$$ c_{x}=c+N_{x}(C) \in C_{x} $$
for the image of $c$ under the quotient projection $C\ra C_{x}$
Note that we have $R_{x}=(\Cs{M})_{x}$ if we regard $R$ as an $R$-module.
Moreover, the action of $R$ on $C_{x}$ factors through
$R_{x}$, so
we can regard $C_{x}$ as an $R_{x}$-module.
If $D$ is another $R$-module and $\alpha:C\ra D$ an $R$-linear map,
then $\alpha(N_{x}(C))\subset N_{x}(D)$ and hence $\alpha$ induces
an $R$-linear (and also $R_{x}$-linear) map
$$ \alpha_{x}:C_{x}\lra D_{x}\;.$$

\begin{prop} \label{prop5}
For any Hausdorff manifold $M$ and for any module $C$ over $R=\Cc(M)$ we have:
\begin{enumerate}
\item [(i)]
      If $c\in C$ satisfies $c_{x}=0$ for any $x\in M$, then $c=0$.
\item [(ii)]
      If $(f_{i})\subset R$ is a sequence satisfying $f_{i}f_{i+1}=f_{i+1}$ and
      $\bigcap_{i}\supp(f_{i})=\{x\}$, then
      $$ N_{x}(C)=\bigcup_{i}(1-f_{i})C\;.$$
\item [(iii)]
      The exact sequence
      $$ 0\lra N_{x}(C)\lra C\lra C_{x}\lra 0 $$
      is locally split, i.e. for any finitely generated $R$-submodule
      $B$ of $N_{x}(C)$ there exists an $R$-linear map
      $C\ra N_{x}(C)$ which is identity on $B$. In particular,
      the sequence
      $$ 0\lra N_{x}(C)\ten_{R}D\lra C\ten_{R}D\lra C_{x}\ten_{R}D\lra 0 $$
      is exact for any $R$-module $D$.
\item [(iv)]
      The homomorphisms $N_{x}(C)\ten_{R} C\ra C\ten_{R}C$ and
      $C\ten_{R} N_{x}(C)\ra C\ten_{R}C$
      are injective, and the kernel of $C\ten_{R}C\ra C_{x}\ten_{R}C_{x}$
      equals to
      $$ N_{x}(C\ten_{R}C)=N_{x}(C)\ten_{R} C=C\ten_{R} N_{x}(C)\;.$$
\end{enumerate}
\end{prop}
\begin{proof}
(i)
Choose  a function $\rho\in\Cc(M)$ with $\rho c=c$.
Furthermore, for any $x\in M$ choose a function $p^{x}\in\Cc(M)$ which
equals $1$ on an open neighbourhood $V_{x}\subset M$ of $x$
such that $p^{x}c=0$. Since $(V_{x})_{x\in M}$ is an open cover of $M$,
we can write $\rho=\sum_{i=1}^{n}\rho_{i}$ for some
$\rho_{i}\in\Cc(V_{x_{i}})$ and $x_{i}\in M$,
$i=1,\ldots,n$. In particular we have $\rho_{i}p^{x_{i}}=\rho_{i}$ and hence
$\rho_{i}c=\rho_{i}p^{x_{i}}c=0$.
This implies that $c=\rho c=0$.

(ii)
Since $f_{i}f_{i+1}=f_{i+1}$, it follows that $f_{i}$ equals $1$ on
$\supp(f_{i+1})$ and hence $\supp(f_{i})$ is a neighbourhood of
$\supp(f_{i+1})$. Now if
$c\in N_{x}(C)$ there exists $f\in R$ with $f_{x}=1$ such that $fc=0$.
In particular,
there exists an open neighbourhood $U$ of $x$ such that $f$ equals $1$ on $U$.
Because $\bigcap_{i}\supp(f_{i})=\{x\}$ we can choose $i\geq 1$ such that
$\supp(f_{i})\subset U$, and therefore $f_{i}f=f_{i}$. Thus
$$ c=c-f_{i}fc=(1-f_{i})c\;.$$
Conversely, if $c=(1-f_{i})d$ for some $d\in D$ and some $i\geq 1$, then
$$ f_{i+1}c=f_{i+1}(1-f_{i})d=f_{i+1}d-f_{i+1}d=0\;,$$
so $c\in N_{x}(C)$.

(iii)
We can choose a sequence $(f_{i})\subset R$ satisfying the condition in (ii).
Then (ii) implies that
for a finitely generated $R$-submodule $B$ of $N_{x}(C)$ there exists $i\geq 1$
such that $B\subset (1-f_{i})C$.
Let $\alpha:C\ra N_{x}(C)$ be given by $\alpha(c)=(1-f_{i+1})c$.
For any $c\in C$ we have
$$ \alpha((1-f_{i})c)=(1-f_{i})(1-f_{i+1})c=(1-f_{i})c\;,$$
so $\alpha$ is identity on $(1-f_{i})C$ and hence also on $B$.

For any  $\sum_{j=1}^{k}c_{j}\ten d_{j}\in\Ker(N_{x}(C)\ten_{R}D\ra C\ten_{R}D)$,
let  $B$ be the $R$-submodule generated by $c_{1},\ldots,c_{k}\in N_{x}(C)$,
and choose $\alpha:C\ra N_{x}(C)$ which is identity on $B$, as above. Then we have
$$ 0=(\alpha\ten\id)(\sum_{j=1}^{k}c_{j}\ten d_{j})
   =\sum_{j=1}^{k}c_{j}\ten d_{j}\;.$$
This shows that $N_{x}(C)\ten_{R}D\ra C\ten_{R}D$ is injective.

(iv)
The first part of the statement is a direct consequence of (iii).
It is then well known that the kernel of
$C\ten_{R}C\ra C_{x}\ten_{R}C_{x}$ equals $N_{x}(C)\ten_{R} C+C\ten_{R} N_{x}(C)$,
so we only need to prove that
$$ N_{x}(C\ten_{R}C)=N_{x}(C)\ten_{R} C=C\ten_{R} N_{x}(C)\;.$$
If $\sum_{j=1}^{k}c_{j}\ten d_{j}\in C\ten_{R}N_{x}(C)$, we can choose
$f\in\Cc(M)$ with
$f_{x}=1$ such that $fd_{j}=0$ for any $1\leq j\leq k$. It follows that
$$ f(\sum_{j=1}^{k}c_{j}\ten d_{j})=0 \;,$$
so $\sum_{j=1}^{k}c_{j}\ten d_{j}\in N_{x}(C\ten_{R}C)$. This shows that
$C\ten_{R} N_{x}(C)\subset N_{x}(C\ten_{R}C)$. Conversely, let
$\sum_{j=1}^{k}c_{j}\ten d_{j}\in N_{x}(C\ten_{R}C)$. Then there exists
$p\in\Cc(M)$ with $p_{x}=1$ such that
$p(\sum_{j=1}^{k}c_{j}\ten d_{j})=\sum_{j=1}^{k}c_{j}\ten pd_{j}=0$,
and therefore
$$ \sum_{j=1}^{k}c_{j}\ten d_{j}=\sum_{j=1}^{k}c_{j}\ten (d_{j}-pd_{j})
   \in C\ten_{R}N_{x}(C) \;.$$
Indeed, we can choose $q\in\Cc(M)$ with $q_{x}=1$ and $qp=q$, so
$q(d_{j}-pd_{j})=0$ and hence $d_{j}-pd_{j}\in N_{x}(C)$. This proves
that $C\ten_{R} N_{x}(C)= N_{x}(C\ten_{R}C)$. The equality
$N_{x}(C)\ten_{R}C = N_{x}(C\ten_{R}C)$ is obtained by an analogous argument.
\end{proof}

\begin{prop}  \label{prop6}
Let $M$ be a Hausdorff manifold, let $(C,\cm,\cu)$ be a coalgebra
over $R=\Cc(M)$ and let $x\in M$. Then
$$ (C_{x},\cm_{x},\cu_{x}) $$
is a coalgebra over the ring $R_{x}=(\Cs{M})_{x}$,
called the local coalgebra of $C$ at $x$.
If $D$ is another coalgebra over $\Cc(M)$ and
if $\alpha:C\ra D$ is a homomorphism of coalgebras, then
$$ \alpha_{x}:C_{x}\lra D_{x} $$
is also a homomorphism of coalgebras.
\end{prop}
\begin{proof}
The counit $\cu$ induces the $R$-linear map
$$ \cu_{x}:C_{x}\lra R_{x}\;,$$
while the comultiplication $\cm$ induces
$$ \cm_{x}:C_{x}\lra (C\ten_{R}C)_{x}=C_{x}\ten_{R_{x}} C_{x}\;.$$
Note that we can here identify
$(C\ten_{R}C)_{x}$ and $C_{x}\ten_{R} C_{x}$ by Proposition \ref{prop5} (iv),
and that $C_{x}\ten_{R} C_{x}=C_{x}\ten_{R_{x}} C_{x}$.
It is now easy to check that $(\cm_{x},\cu_{x})$ is
a coalgebra structure on $C_{x}$ over $R_{x}$. Finally,
the equations
$$ \cu_{x}\com\alpha_{x}=(\cu\com\alpha)_{x}=\cu_{x} $$
and
$$ \cm_{x}\com\alpha_{x}=(\cm\com\alpha)_{x}
   =((\alpha\ten\alpha)\com\cm)_{x}=(\alpha\ten\alpha)_{x}\com\cm_{x}
   =(\alpha_{x}\ten\alpha_{x})\com\cm_{x} $$
show that $\alpha_{x}$ is a homomorphisms of coalgebras.
\end{proof}

Let $C$ be a coalgebra over $R=\Cc(M)$.
We shall say that an element $c$ of $C$ is {\em normalized at} $x\in M$ if
$\cu(c)_{x}=1$. The element $c\in C$
is {\em normalized on} a subset $W$ of $M$ if it
is normalized at any point of $W$.

Since the ring $R_{x}=(\Cs{M})_{x}$ has the identity, we can consider the set
$G(C_{x})$ of {\em grouplike} elements of $C_{x}$,
$$ G(C_{x})=\{\zeta\in C_{x}\,|\, \cm_{x}(\zeta)
   =\zeta\ten\zeta\,,\;\cu_{x}(\zeta)=1\}\;.$$
The following proposition compares grouplike elements of $C_{x}$
with the weakly grouplike
elements of $C$ defined in Section \ref{section:sc}:

\begin{prop}  \label{prop7}
Let $M$ be a Hausdorff manifold, let $C$ be a coalgebra over $R=\Cc(M)$, and
let $x\in M$. Then we have:
\begin{enumerate}
\item [(i)]
      If $c\in G_{w}(C)$ is normalized at $x$, then $c_{x}\in G(C_{x})$.
\item [(ii)]
      If $c,d\in G_{w}(C)$ are normalized on open neighbourhoods $U$
      respectively
      $V$ of $x$ and if $c_{x}=d_{x}$, then there exists an open neighbourhood
      $W\subset U\cap V$ of $x$ such that $c_{y}=d_{y}$ for any $y\in W$.
\item[(iii)]
      If $c$ is an element of $C$ such that $c_{x}\in G(C_{x})$, then
      $c$ is normalized on an open neighbourhood $U$ of $x$, and
      there exists $f\in\Cc(U)$ with $f_{x}=1$ such that $fc\in G_{w}(C)$.
\end{enumerate}
\end{prop}
\begin{proof}
(i)
Let $\cm(c)=c\ten c'$ for some $c'\in C$. We have $\cu_{x}(c_{x})=\cu(c)_{x}=1$
and $\cm_{x}(c_{x})=c_{x}\ten c'_{x}$. The counit property implies that
$c_{x}=\cu_{x}(c_{x})c'_{x}=c'_{x}$, so $c_{x}\in G(C_{x})$.

(ii)
Since $c_{x}=d_{x}$ there exists $f\in\Cc(M)$ with $f_{x}=1$ such that $fc=fd$.
Choose an open neighbourhood $W\subset U\cap V$ of $x$ so small that
$f$ equals $1$ on $W$. Then for any $y\in W$ we have
$c_{y}=(fc)_{y}=(fd)_{y}=d_{y}$.

(iii)
Note that $\cu(c)_{x}=\cu_{x}(c_{x})=1$, so the set $U$ of all $y\in M$ with
$\cu(c)_{y}=1$ is an open neighbourhood of $x$.
Since $(c\ten c-\cm(c))_{x}=0$,
Proposition \ref{prop5} (iv) implies that
$$ c\ten c-\cm(c)=\sum_{i=1}^{n}c_{i}\ten d_{i} $$
for some $c_{i}\in C$ and $d_{i}\in N_{x}(C)$, $i=1,\ldots,n$.
We can now choose $f\in\Cc(U)$ with $f_{x}=1$ such that
$fd_{i}=0$ for all $1\leq i\leq n$. This yields
$f(c\ten c-\cm(c))=0$
and hence $\cm(fc)=fc\ten c$, so $fc\in G_{w}(C)$.
\end{proof}

Let $C$ be a coalgebra over $R=\Cc(M)$.  We define
the {\em spectral sheaf} $\Esp(C)$
of the coalgebra $C$ to be the
sheaf $\pi_{sp}(C):\Esp(C)\ra M$ on $M$ with the stalk
$$ \Esp(C)_{x}=G(C_{x}) $$
at any $x\in M$.
The topology on $\Esp(C)$ is given by a basis consisting
of $\pi_{sp}(C)$-elementary subsets of $\Esp(C)$ of the form
$$ c_{W}=\{ c_{x}\in G(C_{x})\,|\, x\in W\}\subset \Esp(C)\;,$$
where $W$ is any open subset of $M$ and
$c\in G_{w}(C)$ is any weakly grouplike element of $C$ normalized
on $W$. Proposition \ref{prop7} implies that this is indeed a basis
for topology on $\Esp(C)$, so $\Esp(C)$ is a sheaf on $M$.

Let $D$ be another coalgebra over $R=\Cc(M)$ and let $\alpha:C\ra D$
be a homomorphism of coalgebras. By Proposition \ref{prop6} the homomorphism
$\alpha$ induces  the homomorphism of coalgebras over $R_{x}=(\Cs{M})_{x}$
$$ \alpha_{x}:C_{x}\lra D_{x}\;, $$
for any $x\in M$. In particular, $\alpha_{x}(G(C_{x}))\subset G(D_{x})$.
We can now define a map over $M$
$$ \Esp(\alpha):\Esp(C)\lra \Esp(D)\;, $$
given on the stalk of $\Esp(C)$ at any point $x$ of $M$ by
$$ \Esp(\alpha)|_{\Esp(C)_{x}}
   =\alpha_{x}|_{\Esp(C)_{x}}:\Esp(C)_{x}\lra \Esp(D)_{x} \;.$$
To see that $\Esp(\alpha)$ is a morphism of sheaves, it is enough
to observe that
for a basic open $c_{W}$, given by an open subset $W$ of $M$ and a weakly
grouplike $c\in G_{w}(C)$ normalized on $W$, the element $\alpha(c)$
is also weakly grouplike and normalized on $W$, and
$$ \Esp(\alpha)(c_{W})=\alpha(c)_{W}\;.$$
We may conclude that we defined a functor
$$ \Esp:\CoAlg(\Cc(M))\lra \Sh(M)\;.$$

\begin{theo}  \label{theo8}
Let $M$ be a Hausdorff manifold and
$\pi:E\ra M$ a sheaf on $M$.
There is a natural isomorphism of sheaves on $M$
$$ \Phi_{E}: E\lra \Esp(\Cc(E)) $$
defined for any $y\in E$ by
$$ \Phi_{E}(y)=(f\com \pi|_{U})_{\pi(y)}\;,$$
where $U$ is any $\pi$-elementary subset  of $E$ with $y\in U$ and $f$
is any function in $\Cc(\pi(U))$ which satisfies $f_{\pi(y)}=1$.
\end{theo}
\begin{proof}
By Proposition \ref{prop4} we know that $f\com\pi|_{U}$
is a weakly grouplike element of $\Cc(E)$. Furthermore, the grouplike element
$(f\com\pi|_{U})_{\pi(y)}$ of $\Cc(E)_{\pi(y)}$
does not depend on the choice of  $U$ and $f$,
so $\Phi_{E}$ is well defined. The map $\Phi_{E}$ is surjective, because
any element of $\Esp(\Cc(E))$ over $x\in M$
equals $u_{x}$ for some $u\in G_{w}(\Cc(E))$ normalized at $x$,
by Proposition \ref{prop7} (iii). Proposition \ref{prop4}
shows that $u\in\Cc(U)$ for some $\pi$-elementary subset $U$
of $E$, so $u=\cu(u)\com\pi|_{U}$. Now
$$ u_{x}=(\cu(u)\com \pi|_{U})_{x}=\Phi_{E}(y)\;,$$
where $y$ is the unique element of $U$ with $\pi(y)=x$.

To see that $\Phi_{E}$ is injective, observe that
if $U'$ is another  $\pi$-elementary
subset of $E$, $f'\in\Cc(\pi(U'))$ and $y'\in U'$ such that
$(f\com\pi|_{U})_{y}=(f'\com\pi|_{U'})_{y'}$,
then $\pi(y)=\pi(y')$ and there exists a function $p\in\Cc(M)$
with $p_{\pi(y)}=1$ such that
$pf\com\pi|_{U}=pf'\com\pi|_{U'}$. In particular
$y\in U'$ and hence $y=y'$.

To prove that $\Phi_{E}$ is a diffeomorphism it
is now enough to show that (small) $\pi$-elementary subsets
of $E$ map to $\pi_{sp}(\Cc(E))$-elementary subsets of the sheaf
$\Esp(\Cc(E))$.
Take any $\pi$-elementary subset $U$ of $E$ and a function
$f\in\Cc(\pi(U))$ which equals $1$ on an open subset $W$ of $\pi(U)$.
Now $\Phi_{E}(U\cap\pi^{-1}(W))=(f\com\pi|_{U})_{W}$.

Finally, we have to check that $\Phi_{E}$ is natural in $E$.
Take any morphism of sheaves $\phi:E\ra F$ over $M$, let $y\in E$ and
choose a $\pi$-elementary subset $U$ of $E$ with $y\in U$
and a function $f\in\Cc(\pi(U))$ with
$f_{\pi(y)}=1$. Note that $\phi(U)$ is $\xi$-elementary, where $\xi$ is
the projection of the sheaf $F$. Now we have
\begin{eqnarray*}
\Esp(\phi_{+})(\Phi_{E}(y))
& = & \Esp(\phi_{+})(f\com\pi|_{U})_{\pi(y)}
      = (\phi_{+})_{\pi(y)}(f\com\pi|_{U})_{\pi(y)} \\
& = & (\phi_+{}(f\com\pi|_{U}))_{\pi(y)}
      = (f\com\xi|_{\phi(U)})_{\pi(y)} \\
& = & \Phi_{F}(\phi(y))\;.
\end{eqnarray*}
This proves that $\Phi:\id\ra\Esp\com\Cc$ is a natural isomorphism.
\end{proof}

\begin{theo}  \label{theo9}
Let $M$ be a Hausdorff manifold and let $C$ be a coalgebra over $\Cc(M)$.
Then there is a natural homomorphism of coalgebras over $\Cc(M)$
$$ \Psi_{C}:\Cc(\Esp(C))\lra C $$
given by
$$ \Psi_{C}(\sum_{i=1}^{k}f_{i}\com \pi_{sp}(C)|_{(c_{i})_{W_{i}}})
   =\sum_{i=1}^{k}f_{i}c_{i} \;,$$
where $c_{i}$ is a weakly grouplike element of $C$ normalized on an open
subset $W_{i}$ of $M$ and $f_{i}\in\Cc(W_{i})$, for any $i=1,\ldots,k$.
\end{theo}
\begin{proof}
To simplify the notation, write $\pi=\pi_{sp}(C)$.
First we have to prove that $\Psi_{C}$ is well defined.
Any function in $\Cc(\Esp(C))$ is a sum of functions
of the form $f\com\pi|_{c_{W}}$, for $c\in G_{w}(C)$ normalized on
an open subset $W$ of $M$ and $f\in\Cc(W)$,
because the $\pi$-elementary subsets
of the form $c_{W}$
constitute an open cover (in fact a basis) of $\Esp(C)$.
Suppose that for any $i=1,\ldots,k$
we have a weakly grouplike
element $c_{i}\in G_{w}(C)$
normalized on an open subset $W_{i}$ of $M$ and a function
$f_{i}\in\Cc(W_{i})$ such that
$$ u=\sum_{i=1}^{k}f_{i}\com \pi|_{(c_{i})_{W_{i}}}=0\;.$$
We have to show that $\Psi_{C}(u)=\sum_{i=1}^{k}f_{i}c_{i}=0$.

By Proposition \ref{prop5} (i) it is sufficient to prove that
$\sum_{i=1}^{k}(f_{i})_{x}(c_{i})_{x}=0$ for any $x\in M$.
Take any $x\in M$ and let $\zeta_{1},\ldots,\zeta_{n}$
be all distinct elements of $G(C_{x})$ in
$\bigcup_{i=1}^{k}(c_{i})_{W_{i}}$.
For any $1\leq l\leq n$
denote by $I_{l}$ the set
of all $1\leq i\leq k$ with $x\in W_{i}$ and
$\zeta_{l}=(c_{i})_{x}$. Let $I_{0}$ be the set of all
$1\leq i\leq k$ with $x\not\in W_{i}$. Note that
$\{1,\ldots,k\}$ is the disjoint union of $I_{0},\ldots,I_{n}$.
Furthermore, we have $(f_{i})_{x}=0$ for any $i\in I_{0}$, while
$\sum_{i\in I_{l}}(f_{i})_{x}=0$
because $u_{\zeta_{l}}=0$, for any $1\leq l\leq n$.
This yields
$$
\sum_{i=1}^{k}(f_{i})_{x}(c_{i})_{x}
= \sum_{i\in I_{0}}(f_{i})_{x}(c_{i})_{x}
  +\sum_{l=1}^{n} \sum_{i\in I_{l}} (f_{i})_{x}(c_{i})_{x}
= \sum_{l=1}^{n} (\sum_{i\in I_{l}} (f_{i})_{x})\zeta_{l} =0\;.$$

Now we have to check that $\Psi_{C}$ is a homomorphism of coalgebras.
Let $W$ be an open subset of $M$, let $c$ be a weakly grouplike
element of $C$ normalized on $W$ and let $f\in\Cc(W)$. Choose also
a function $\eta\in\Cc(W)$
which satisfies $\eta f=f$. Since $\cm(c)=c\ten c'=c'\ten c$ for
some $c'\in C$,
it follows that $\cm(fc)=fc'\ten c=f\cu(c)c'\ten c=fc\ten c$.
Therefore we have
\begin{eqnarray*}
\lefteqn{   (\Psi_{C}\ten\Psi_{C})(\cm(f\com \pi|_{c_{W}}))    }
&   & \\
& = & (\Psi_{C}\ten\Psi_{C})(f\com \pi|_{c_{W}}\ten \eta\com \pi|_{c_{W}}) \\
& = & fc\ten\eta c=\eta fc\ten c=fc\ten c = \cm(fc) \\
& = & \cm(\Psi_{C}(f\com \pi|_{c_{W}}))\;.
\end{eqnarray*}
Furthermore, we have
$$ \cu(\Psi_{C}(f\com \pi|_{c_{W}}))=\cu(fc)=f=\cu(f\com \pi|_{c_{W}})\;.$$

Finally, let $\alpha:C\ra D$ be a homomorphism of coalgebras over $\Cc(M)$.
Note that $\alpha(c)$ is weakly grouplike and normalized on $W$. We have
\begin{eqnarray*}
\Psi_{D}(\Esp(\alpha)_{+}(f\com\pi|_{c_{W}}))
& = & \Psi_{D}(f\com\pi_{sp}(D)|_{\Esp(\alpha)(c_{W})}) \\
& = & \Psi_{D}(f\com\pi_{sp}(D)|_{\alpha(c)_{W}}) \\
& = & f\alpha(c)=\alpha(fc) \\
& = & \alpha(\Psi_{C}(f\com\pi|_{c_{W}}))\;,
\end{eqnarray*}
so $\Psi:\Cc\com\Esp\ra\id$ is a natural transformation.
\end{proof}

\begin{theo}  \label{theo10}
For any Hausdorff manifold $M$,
the functor 
$$ \Cc:\Sh(M)\lra\CoAlg(M) $$
is fully faithful and left
adjoint to the functor $\Esp$.
\end{theo}
\begin{proof}
The unit of the adjunction is $\Phi$ of Theorem \ref{theo8},
while $\Psi$ of Theorem \ref{theo9}
is the counit. Indeed, take any sheaf $\pi:E\ra M$,
a $\pi$-elementary subset $U$ of $E$ and any $f\in\Cc(\pi(U))$.
Choose also a function $\eta\in\Cc(\pi(U))$ which equals $1$
on an open neighbourhood $W\subset U$ of $\supp(f)$, and put
$V=(\pi|_{U})^{-1}(W)$. Now
\begin{eqnarray*}
\Psi_{\Cc(E)}((\Phi_{E})_{+}(f\com\pi|_{U}))
& = & \Psi_{\Cc(E)}((\Phi_{E})_{+}(f\com\pi|_{V})) \\
& = & \Psi_{\Cc(E)}(f\com\pi_{sp}(\Cc(E))|_{\Phi_{E}(V)}) \\
& = & \Psi_{\Cc(E)}(f\com\pi_{sp}(\Cc(E))|_{(\eta\com\pi|_{U})_{W}}) \\
& = & f(\eta\com\pi|_{U})=f\com\pi|_{U}\;.
\end{eqnarray*}
Next, take a coalgebra $C$ over $\Cc(M)$ and any $c_{x}\in\Esp(C)_{x}$
represented
by a weakly grouplike $c\in C$ normalized on an open neighbourhood $W$
of $x\in M$.
Choose also a function $f\in\Cc(W)$ with $f_{x}=1$. We have
\begin{eqnarray*}
\Esp(\Psi_{C})(\Phi_{\Esp(C)}(c_{x}))
& = & \Esp(\Psi_{C})((f\com\pi_{sp}(C)|_{c_{W}})_{x}) \\
& = & (\Psi_{C}(f\com\pi_{sp}(C)|_{c_{W}}))_{x} \\
& = & (fc)_{x}=c_{x}\;.
\end{eqnarray*}
This proves that $\Esp$ is right adjoint to $\Cc$.
Since the unit of this adjunction
is a natural isomorphism by Theorem \ref{theo8}, it follows that
$\Cc$ is fully faithful.
\end{proof}

Let $M$ be a Hausdorff manifold and let $C$ be a coalgebra over $\Cc(M)$.
We say that a subset $\Lambda\subset C$ {\em normally generate} $C$
if any $c\in C$ can be written as a sum
$$ c=\sum_{i=1}^{k}f_{i}c_{i} $$
for some $f_{1},\ldots,f_{k}\in\Cc(M)$ and $c_{1},\ldots,c_{k}\in \Lambda$
such that $c_{i}$ is normalized on $\supp(f_{i})$ for any $i=1,\ldots,k$.

We say that $\Lambda\subset C$
is {\em normally linearly independent at} $x\in M$
if for any $c_{1},\ldots,c_{k}\in\Lambda$ normalized at $x$
and for any $f_{1},\ldots,f_{k}\in\Cc(M)$ satisfying
$$ \sum_{i=1}^{k}f_{i}c_{i}=0 \;,$$
it follows that either $(f_{i})_{x}=0$ for any $i=1,\ldots,k$ or
there exist distinct $1\leq i,j\leq k$ and a function $f\in\Cc(M)$ with
$f_{x}=1$ such that $fc_{i}=fc_{j}$.
A subset $\Lambda$ is {\em normally linearly independent} if it is
normally linearly independent at any $x\in M$.

For any subset $\Lambda\subset C$ and any $x\in M$ let us denote
$$ \Lambda^{\cu}_{x}
   =\{c_{x}\,|\,c\in\Lambda\,,\;\cu(c)_{x}=1\} \subset C_{x}\;.$$

\begin{prop}  \label{prop11}
Let $M$ be a Hausdorff manifold and let $C$ be a coalgebra over $\Cc(M)$.
\begin{enumerate}
\item [(i)]
      A subset $\Lambda\subset C$ normally generate $C$ if and only if
      the set $\Lambda^{\cu}_{x}$
      generate $C_{x}$ as an $(\Cs{M})_{x}$-module, for any $x\in M$.
\item [(ii)]
      A subset $\Lambda\subset C$ is normally linearly independent at
      $x\in M$ if the set $\Lambda^{\cu}_{x}$
      is $(\Cs{M})_{x}$-linearly independent in $C_{x}$.
\end{enumerate}
\end{prop}
\begin{proof}
(i)
Suppose that $\Lambda$ normally generate $C$.
Take any $x\in M$ and $c_{x}\in C_{x}$ represented by an element
$c\in C$, and write
$c=\sum_{i=1}^{k}f_{i}c_{i}$
for some $f_{1},\ldots,f_{k}\in\Cc(M)$ and some $c_{1},\ldots,c_{k}\in\Lambda$
such that $c_{i}$ is normalized on $\supp(f_{i})$ for any $i=1,\ldots,k$.
Let $I$ be the set of those $1\leq i\leq k$ for which
$x\in\supp(f_{i})$. Now we have
$$ c_{x}=\sum_{i\in I}(f_{i})_{x}(c_{i})_{x}\;,$$
with $\cu(c_{i})_{x}=1$ for any $i\in I$.

Conversely, assume that $\Lambda^{\cu}_{x}$ generate $C_{x}$
for any $x\in M$, and take any $c\in C$. For any $x\in M$
we can write
$$ c_{x}=\sum_{i=1}^{k_{x}}(f^{x}_{i})_{x}(c^{x}_{i})_{x} $$
for some functions $f^{x}_{i}\in\Cc(M)$ and
some $(c^{x}_{i})_{x}\in \Lambda^{\cu}_{x}$ represented by
$c^{x}_{i}\in\Lambda$ normalized at $x$,
for any $i=1,\ldots,k_{x}$.
Therefore there exists $f^{x}\in\Cc(M)$ which equals $1$ on an
open neighbourhood $V_{x}$ of $x$ such that
$$ f^{x}c=\sum_{i=1}^{k_{x}}f^{x}f^{x}_{i}c^{x}_{i}\;.$$
We can choose $V_{x}$ so small that each $c^{x}_{i}$ is normalized
on $V_{x}$.

Choose $\rho\in\Cc(M)$ with $\rho c=c$.
Since $(V_{x})_{x\in M}$ is an open cover of $M$, we can write
$\rho=\sum_{j=1}^{m}\rho_{j}$ for some $\rho_{j}\in\Cc(V_{x_{j}})$ and
$x_{j}\in M$, $j=1,\ldots,m$. In particular we have
$\rho_{j}f^{x_{j}}=\rho_{j}$. It follows that
$$ c=\rho c=\sum_{j=1}^{m}\rho_{j}f^{x_{j}}c
   =\sum_{j=1}^{m}\rho_{j}\sum_{i=1}^{k_{x_{j}}}f^{x_{j}}f^{x_{j}}_{i}c^{x_{j}}_{i}\;.$$
Since each $c^{x_{j}}_{i}\in\Lambda$ is normalized on $\supp(\rho_{j})$ and hence
also on $\supp(\rho_{j}f^{x_{j}}f^{x_{j}}_{i})$, it follows that
$\Lambda$ normally generate $C$.

(ii)
Suppose first that $\Lambda^{\cu}_{x}$ is $(\Cs{M})_{x}$-linearly independent.
Let $c_{1},\ldots,c_{k}$ be elements of $\Lambda$
normalized at $x$ and let $f_{1},\ldots,f_{k}\in\Cc(M)$ be such that
$$ \sum_{i=1}^{k}f_{i}c_{i}=0 \;.$$
If all $(c_{i})_{x}$ are distinct then it follows that
$(f_{i})_{x}=0$ for any $i=1,\ldots,k$. On the other hand, if there
exist distinct $1\leq i,j\leq k$ with $(c_{i})_{x}=(c_{j})_{x}$ then by definition
there is a function $f\in\Cc(M)$ with $f_{x}=1$ such that $fc_{i}=fc_{j}$.

Conversely, assume that $\Lambda$ is normally linearly independent at $x$.
Take any distinct $(c_{1})_{x},\ldots,(c_{k})_{x}\in\Lambda^{\cu}_{x}$
represented by $c_{1},\ldots,c_{k}\in\Lambda$ normalized at $x$.
Let $f_{1},\ldots,f_{k}\in \Cc(M)$ be such that
$$ \sum_{i=1}^{k}(f_{i})_{x}(c_{i})_{x}=0\;.$$
Hence there exists $f\in\Cc(M)$ with $f_{x}=1$ such that
$$ \sum_{i=1}^{k}ff_{i}c_{i}=0\;.$$
Since $(c_{i})_{x}$ are distinct, it follows that
$(f_{i})_{x}=(ff_{i})_{x}=0$ for any $i=1,\ldots,k$.
\end{proof}

Let $M$ be a Hausdorff manifold, let $C$ be a coalgebra over $\Cc(M)$, and
let $\Lambda\subset G_{w}(C)$. Define a subsheaf $\EspL(C)$ of the spectral
sheaf $\Esp(C)$ on $M$ by
$$ \EspL(C)_{x}=\{c_{x}\,|\,c\in\Lambda\,,\;\cu(c)_{x}=1\}=\Lambda^{\cu}_{x} $$
for any $x\in M$. This is indeed an open subspace of $\Esp(C)$ because
for each $c_{x}\in\EspL(C)_{x}$ represented by $c\in\Lambda$ with $\cu(c)_{x}=1$ there
exists an open neighbourhood $W$ of $x$ such that $c$ is normalized on $W$, and hence
$c_{W}\subset \EspL(C)$. Denote by $\inc:\EspL(C)\ra\Esp(C)$ the inclusion.
The induced map $\inc_{+}:\Cc(\EspL(C))\ra\Cc(\Esp(C))$ is a monomorphism given by
extension by zero. Let us denote by
$$ \Psi_{C}^{\Lambda}=\Psi_{C}\com\inc_{+}:\Cc(\EspL(C))\lra C $$
the restriction of $\Psi_{C}$ to $\Cc(\EspL(C))$.

\begin{theo}  \label{theo12}
Let $M$ be a Hausdorff manifold, let $C$ be a coalgebra
over $\Cc(M)$ and let $\Lambda$ be subset of $G_{w}(C)$. Then we have:
\begin{enumerate}
\item [(i)]
      $\Psi_{C}^{\Lambda}$ is an epimorphism if and only if
      $\Lambda$ normally generate $C$.
\item [(ii)]
      $\Psi_{C}^{\Lambda}$ is a monomorphism if and only if
      $\Lambda$ is normally linearly independent.
\end{enumerate}
\end{theo}
\Rem
In particular, this theorem implies that $\Psi_{C}$ is an epimorphism
if and only if $G_{w}(C)$ normally generate $C$, and that
$\Psi_{C}$ is a monomorphism if and only if $G_{w}(C)$ is
normally linearly independent.

\begin{proof}
Write $\pi=\pi_{sp}(C)$.
First note that the $\pi$-elementary subsets $c_{W}$,
for $c\in\Lambda$ normalized on an open subset $W\subset M$,
form a basis on $\EspL(C)$.
Therefore any element of $\Cc(\EspL(C)$ can be written as a sum of
functions of the form $f\com \pi|_{c_{W}}$, where
$c\in\Lambda$ is normalized on $W$ and $f\in\Cc(W)$.
Since  $\Psi_{C}^{\Lambda}$ is given by
$$ \Psi_{C}^{\Lambda}(f\com \pi|_{c_{W}})=fc\;,$$
this clearly implies (i).

To prove (ii), suppose first that
$\Psi_{C}^{\Lambda}$ is injective.  Take any
$c_{1},\ldots,c_{k}\in \Lambda$ normalized at a point $x\in M$, and
any $f_{1},\dots,f_{k}\in\Cc(M)$ such that $\sum_{i=1}^{k}f_{i}c_{i}=0$.
Assume that all $(c_{1})_{x},\ldots,(c_{k})_{x}$ are distinct.
Let $W$ be an open neighbourhood of $x$ such that each $c_{i}$ is normalized on $W$,
and choose any $\nu\in\Cc(W)$ with $\nu_{x}=1$.
Since
$$ \Psi_{C}^{\Lambda}(\sum_{i=1}^{k}\nu f_{i}\com\pi|_{(c_{i})_{W}})
   =\nu\sum_{i=1}^{k}f_{i}c_{i}=0\;,$$
it follows by injectivity of $\Psi_{C}^{\Lambda}$ that
$$ \sum_{i=1}^{k}\nu f_{i}\com\pi|_{(c_{i})_{W}}=0\;.$$
Thus for any $1\leq j\leq k$ we have
$$ 0=(\sum_{i=1}^{k}\nu f_{i}\com\pi|_{(c_{i})_{W}})_{(c_{j})_{x}}
    =\pi^{\ast}_{(c_{j})_{x}}((\nu f_{j})_{x})\;,$$
and hence $(f_{j})_{x}=(\nu f_{j})_{x}=0$.

Conversely, assume that $\Lambda$ is normally linearly independent at any $x\in M$.
Take any $u\in\Ker(\Psi_{C}^{\Lambda})$ and write
$$ u=\sum_{i=1}^{k} f_{i}\com\pi|_{(c_{i})_{W_{i}}}\;, $$
where each $c_{i}\in\Lambda$ is normalized on an open
set $W_{i}\subset M$ and $f_{i}\in\Cc(W_{i})$.
Thus we have $\Psi_{C}^{\Lambda}(u)=\sum_{i=1}^{k}f_{i}c_{i}=0$.

We have to show that $u=0$. It is sufficient to show that
$u_{(c_{j})_{x}}=0$ for any $x\in M$ and any $1\leq j\leq k$ with $x\in W_{j}$.
To this end, take any $x\in M$ and denote by
$\zeta_{1},\ldots,\zeta_{n}$ all distinct
elements of $\EspL(C)_{x}$ in $\bigcup_{i=1}^{k}(c_{i})_{W_{i}}$.
Let $I_{0}$ be the set of all $1\leq i\leq k$ for which
$x\not\in W_{i}$. Moreover, for any $1\leq l\leq n$ let $I_{l}$ be the set of
all $1\leq i\leq k$ with $\zeta_{l}=(c_{i})_{x}$.
Since $(f_{i})_{x}=0$ for any $i\in I_{0}$, we have
$$ 0=\sum_{i=1}^{k}(f_{i})_{x}(c_{i})_{x}=\sum_{i\in I_{0}}(f_{i})_{x}(c_{i})_{x}
     +\sum_{l=1}^{n}\sum_{i\in I_{l}}(f_{i})_{x}(c_{i})_{x}
     =\sum_{l=1}^{n}(\sum_{i\in I_{l}}(f_{i})_{x})\zeta_{l}\;.$$
Since $\zeta_{1},\ldots,\zeta_{n}$ are
$(\Cs{M})_{x}$-linearly independent
by Proposition \ref{prop11} (ii), it follows that
$$ \sum_{i\in I_{l}}(f_{i})_{x}=0 $$
for any $1\leq l\leq n$. Now for any $1\leq m\leq n$ and for any
$j\in I_{m}$ we have
\begin{eqnarray*}
u_{(c_{j})_{x}}
& = & \sum_{i=1}^{k} (f_{i}\com\pi|_{(c_{i})_{W_{i}}})_{(c_{j})_{x}}
      = \sum_{l=1}^{n}\sum_{i\in I_{l}} (f_{i}\com\pi|_{(c_{i})_{W_{i}}})_{\zeta_{m}} \\
& = & \sum_{i\in I_{m}} \pi^{\ast}_{\zeta_{m}}((f_{i})_{x})
      = \pi^{\ast}_{\zeta_{m}}(\sum_{i\in I_{m}} (f_{i})_{x})=0\;.
\end{eqnarray*}
This proves that $\Psi_{C}^{\Lambda}$ is a monomorphism.
\end{proof}

\begin{cor}  \label{cor13}
Let $M$ be a Hausdorff manifold and let $C$ be a coalgebra over $\Cc(M)$.
If for any $x\in M$ the $(\Cs{M})_{x}$-module
$C_{x}$ is flat and if
any two grouplike elements of $C_{x}$ are
$(\Cs{M})_{x}$-linearly independent, then
$\Psi_{C}:\Cc(\Esp(C))\ra C$ is a monomorphism.
\end{cor}
\begin{proof}
This follows from Theorem \ref{theo12} and Proposition \ref{prop11}
with $\Lambda=G_{w}(C)$, and from
Proposition \ref{prop1} with $L=G(C_{x})$.
\end{proof}

\begin{theo}  \label{theo14}
Let $M$ be a Hausdorff manifold and let $C$ be a coalgebra over $\Cc(M)$.
The following conditions are equivalent:
\begin{enumerate}
\item [(i)]   The coalgebra $C$ is isomorphic to the sheaf coalgebra $\Cc(E)$
              of a sheaf $E$ on $M$.
\item [(ii)]  The homomorphism of coalgebras $\Psi_{C}:\Cc(\Esp(C))\ra C$ is
              an isomorphism.
\item [(iii)] The weakly grouplike elements of $C$ normally generate
              $C$ and are normally linearly independent.
\item [(iv)]  The $(\Cs{M})_{x}$-module $C_{x}$ is freely generated by
              the grouplike elements of $C_{x}$, for any $x\in M$.
\end{enumerate}
Moreover, if the coalgebra $C$ is isomorphic to the sheaf coalgebra $\Cc(E)$
of a sheaf $E$ on $M$, then the sheaf $E$ is isomorphic to
the spectral sheaf $\Esp(C)$ of $C$.
\end{theo}
\begin{proof}
The equivalence between (ii), (iii) and (iv) follows directly from
Theorem \ref{theo12} and Proposition \ref{prop11} with $\Lambda=G_{w}(C)$.
The condition (ii) clearly implies (i),
while (i) implies (ii) by Theorem \ref{theo10}.
The isomorphism between $E$ and $\Esp(C)$ is also a consequence of
Theorem \ref{theo10}.
\end{proof}

\section{Groupoid Hopf algebroids} \label{section:gha}

In this section we show that the space of
smooth functions with compact support
on a non-necessarily Hausdorff
\'{e}tale Lie groupoid over a manifold $M$ admits a natural structure of
a Hopf algebroid over $\Cc(M)$.

Suppose that $A$ is an algebra with a commutative subalgebra $R\subset A$.
Here, and throughout this
paper we assume that algebras are associative,
perhaps without identity, and over a fixed base field $\FF\in\{\RR,\CC\}$.
Assume that $A$ has \textit{local identities in}
$R$, i.e. for any elements $a_{1},\ldots,a_{k}$ of $A$
there exists $r\in R$ such that
$$ a_{i}r=ra_{i}=a_{i}\;, \;\;\;\;\;\;\;\;\;\; i=1,\ldots,k\;.$$
In particular, this implies that $R$ has local identities.
Note that $A$ can be regarded as an $R$-$R$-bimodule, and both
of these module structures are locally unitary.

We shall denote by $A\tenll A$ the tensor product with respect to the left $R$-action
of $R$ on both factors. In other words, $A\tenll A$ is the quotient of
$A\ten A=A\ten_{\FF}A$ over the subspace generated by elements of the form
$ra\ten b-a\ten rb$, for all $a,b\in A$ and $r\in R$.
Analogously, $A\tenrl A$ will denote the tensor product with respect to the
right $R$-action on the first factor and the left $R$-action on the second one.
Note that $A\tenll A$ has a left $R$-action (because $R$ is commutative)
and two natural right $R$-actions, and all these actions commute with each other.

\begin{dfn} \rm \label{dfn15}
Let $A$ be an algebra with local identities in
a commutative subalgebra $R\subset A$.
A \textit{Hopf algebroid structure on} $A$ {\em over} $R$
consists of three $\FF$-linear maps,
the {\em comultiplication}
$$ \cm:A\lra A\tenll A\;,$$
the {\em counit}
$$ \cu:A\lra R $$
and the {\em antipode}
$$ S:A\lra A \;,$$
such that
\begin{enumerate}
\item [(i)]
      $\cm$ and $\cu$ are homomorphisms of left $R$-modules satisfying
      $(\id\ten\cu)\com\cm=\id$ and $(\cu\ten\id)\com\cm=\id$,
\item [(ii)]
      $\cu|_{R}=\id$,
      $\cm|_{R}$ is the canonical embedding $R\cong R\ten_{R} R\subset A\tenll A$,
      and the  two right $R$-actions on $A\tenll A$ coincide on $\cm A$,
\item [(iii)]
      $\cm(ab)=\cm(a)\cm(b)$ for any $a,b\in A$,
\item [(iv)]
      $S|_{R}=\id$ and $S\com S=\id$,
\item [(v)]
      $S(ab)=S(b)S(a)$ for any $a,b\in A$, and
\item [(vi)]
      $\mu\com (S\ten\id)\com\cm=\cu\com S$, where $\mu:A\tenrl A\ra A$ denotes the
      multiplication.
\end{enumerate}
If $R$ is a commutative algebra with local identities, then
a {\em Hopf algebroid over} $R$ is a quadruple $(A,\cm,\cu, S)$, where
$A$ is an algebra which has $R$ for a subalgebra and has local identities in $R$,
and where $(\cm,\cu, S)$ is a Hopf algebroid structure on $A$ over $R$.
\end{dfn}

\Rem
(1) Condition (i) means precisely
that $(A,\cm,\cu)$ is a coalgebra over $R$ with
respect to the left action of $R$ on $A$. A Hopf algebroid
$A=(A,\cm,\cu, S)$ over $R$ is {\em coassociative} respectively
{\em cocommutative} if the coalgebra $(A,\cm,\cu)$
is coassociative respectively cocommutative.

(2)
Note that (ii) implies that the product on $A\ten A$ induces a product
$$ \cm A\ten (A\tenll A)\lra A\tenll A $$
which is used in (iii). The property (iii) together with (ii) yields that
$\cm$ is a homomorphism of right $R$-modules, for both right
$R$-actions on $A\tenll A$. This in particular implies that
one can define a linear map
$\overline{\cm}:A\tenrl A\ra A\tenll A$ with
$$ \overline{\cm}(a\ten b)=\cm(a)(r\ten b)\;,$$
where $r$ is any element of $R$ satisfying $rb=b$.

(3)
In various definitions of Hopf algebroids (or quantum groupoids)
given in the literature \cite{Lu,Mal,Mrc01,Xu},
the base algebra $R$ may be non-commutative and
two maps, the source and the target, from $R$ to $A$ are parts of the structure.
In our definition, we concentrate on a special class of Hopf algebroids
for which $R$ is a commutative subalgebra of $A$ which may not lie in
the center of $A$, while the source and the target maps are both equal to
the inclusion. On the other hand, we do not assume that $A$ has
the identity.

(4)
We do not assume here any direct relation between
$\cm\com S$ and $(S\ten S)\com\cm$ as this is not necessary for our results,
but one should note that these two maps can not be directly compared because
they have different targets. Indeed, the map $ S\ten S$ is an isomorphism
between $A\tenll A$ and $A\tenrr A$. Denote by $\varpi:A\ten A\ra A\tenll A\,\tenrr$
the quotient projection, where $A\tenll A\,\tenrr$ is the tensor product with respect
to both the left and the right actions, i.e. the kernel of $\varpi$ is the subspace
generated by elements of the form $rar'\ten b-a\ten rbr'$, for all $a,b\in A$
and $r,r'\in R$. The definition of an \'{e}tale Hopf algebroid given in \cite{Mrc01}
assumes that
$$ \varpi\com\cm\com S=\varpi\com (S\ten S)\com\cm\;.$$
On the other hand, the comultiplication $\cm$ is injective and $S$ is invertible,
so there exists a unique $\FF$-linear isomorphism $\gamma:\cm A\ra (S\ten S)\cm A$
such that
$$ \gamma\com\cm\com S=(S\ten S)\com\cm\;.$$
This equation, together with the property (vi), implies that
$$ \mu\com (\id\ten S)\com\gamma\com\cm=\cu\;.$$
\vspace{2mm}

Let $R$ be a commutative algebra with local identities.
If $A$ and $B$ are two Hopf algebroids over $R$, then a map
$\alpha:A\ra B$ is a homomorphism of Hopf algebroids over $R$ if $\alpha$ is a
homomorphism of algebras, $\alpha|_{R}=\id$, $\cu\com\alpha=\cu$,
$(\alpha\ten\alpha)\com\cm=\cm\com\alpha$
and $S\com\alpha=\alpha\com S$.
In particular, a homomorphism of Hopf algebroids over $R$ is also a homomorphism
of underlying
coalgebras over $R$. We shall denote by
$$ \HoAlgd(R) $$
the category of Hopf algebroids over $R$ and homomorphisms between them.

Let $A$ be a  Hopf algebroid over $R$. Since $A$ is also a coalgebra
over $R$, we can consider the set of weakly grouplike elements $G_{w}(A)$ of $A$.
A weakly grouplike element $a\in G_{w}(A)$ of $A$
is {\em $S$-invariant} if there exists $a'\in A$ such that
$\cm(a)=a\ten a'$
and
$\cm(S(a))=S(a')\ten S(a)$.
We shall denote the set of $S$-invariant weakly
grouplike elements of $A$ by
$$ G^{S}_{w}(A)\;.$$
The following properties are easy consequences of Definition \ref{dfn15}:

\begin{prop} \label{prop16}
Let $R$ be a commutative algebra with local identities
and let $A$ be a Hopf algebroid over $R$.
Then $R\subset G^{S}_{w}(A)\subset G_{w}(A)$, and both
$G_{w}(A)$ and $G^{S}_{w}(A)$ are closed under the
multiplication in $A$. Furthermore, we have
$S(G^{S}_{w}(A))=G^{S}_{w}(A)$.
If $a\in G^{S}_{w}(A)$ with $\cm(a)=a\ten a'$ and
$\cm(S(a))=S(a')\ten S(a)$, then we have
\begin{enumerate}
\item [(i)]
      $a=\cu(a)a'=\cu(a')a=a\cu(S(a'))=a'\cu(S(a))$,
\item [(ii)]
      $S(a)=\cu(S(a'))S(a)=\cu(S(a))S(a')=S(a')\cu(a)=S(a)\cu(a')$,
\item [(iii)]
      $\cm(a)=\cu(a)a'\ten a'=a'\ten a$,
\item [(iv)]
      $\cm(S(a))=S(a')\ten S(a')\cu(a)=S(a)\ten S(a')$,
\item [(v)]
      $\cu(a)=aS(a')=a'S(a)$ and
\item [(vi)]
      $\cu(S(a))=S(a')a=S(a)a'$.
\end{enumerate}
\end{prop}

\begin{prop} \label{prop17}
Let $R$ be a commutative algebra with local identities
and let $A$ be a Hopf algebroid over $R$. If $A$ is generated
by $G^{S}_{w}(A)$ as an algebra, then
\begin{enumerate}
\item [(i)]
      $A$ is generated by $G^{S}_{w}(A)$ as an abelian group,
\item [(ii)]
      $A$ is cocommutative and coassociative,
\item [(iii)]
      $\cu(ab)=\cu(a\cu(b))$ for any $a,b\in A$,
\item [(iv)]
      $(S\ten\id)\com\overline{\cm}\com (S\ten\id)\com \overline{\cm}=\id$,
\item [(v)]
      $\varpi\com\cm\com S=\varpi\com (S\ten S)\com\cm$, and
\item [(vi)]
      $\cm A$ is generated, as an abelian group,
      by elements of the form $a\ten a'$ satisfying
      $\gamma(a\ten a')=a\ten a'=a'\ten a$.
\end{enumerate}
\end{prop}
\Rem
Here $\overline{\cm}$, $\varpi$ and $\gamma$ are given
as in Remark (2) and (4) to Definition \ref{dfn15}.
This proposition shows that a Hopf algebroid $A$ generated by $G^{S}_{w}(A)$
is an \'{e}tale Hopf algebroid as defined in \cite{Mrc01}.
Note, however, the difference in notation:
in \cite{Mrc01} the coalgebra structure on $A$
is assumed to be over the right $R$-module structure.

\begin{proof}
(i) and (ii) are consequences of Proposition \ref{prop16}.
In fact, recall from Section \ref{section:sc}
that (ii) is true for any coalgebra which is generated, as an abelian group,
by weakly grouplike elements.

(iii)
It is sufficient to prove this for
$a,b\in G^{S}_{w}(A)$. Let $a',b'\in A$ be such that
$\cm(a)=a\ten a'$, $\cm(S(a))=S(a')\ten S(a)$,
$\cm(b)=b\ten b'$ and  $\cm(S(b))=S(b')\ten S(b)$.
Then we have $ab\in G^{S}_{w}$ with
$\cm(ab)=ab\ten a'b'$ and $\cm(S(ab))=S(a'b')\ten S(ab)$,
and also $a\cu(b)\in G^{S}_{w}(A)$ with $\cm(a\cu(b))=a\cu(b)\ten a'$ and
$\cm(S(a\cu(b)))=S(a')\ten S(a\cu(b))$. Therefore
$$ \cu(ab)=abS(a'b')=abS(b')S(a')=a\cu(b)S(a')=\cu(a\cu(b))\;.$$

(iv)
With $a,a'$ as in (iii) and arbitrary $b\in A$ we have
\begin{eqnarray*}
(S\ten\id)\com\overline{\cm}\com (S\ten\id)\com \overline{\cm}(a\ten b)
& = & (S\ten\id)\com\overline{\cm}\com (S\ten\id)(a\ten a'b) \\
& = & (S\ten\id)\com\overline{\cm}(S(a)\ten a'b) \\
& = & (S\ten\id)(S(a')\ten S(a)a'b) \\
& = & a'\ten \cu(S(a))b=a'\cu(S(a))\ten b \\
& = & a\ten b\;.
\end{eqnarray*}

(v) and (vi) are obvious.
\end{proof}

We shall now describe our main example of a Hopf algebroid.
In \cite{Mrc01} we showed that the Connes convolution algebra
of a Hausdorff \'{e}tale Lie groupoid (\cite{Con78,Con82}) has a
natural structure of a cocommutative coassociative Hopf algebroid.
We shall now generalize this construction to non-Hausdorff \'{e}tale
Lie groupoids.

Let $M$ be a Hausdorff manifold.
Recall that an \'{e}tale Lie groupoid over $M$
(see e.g. \cite{CraMoe,Est,Hae01,Mrc99})
is a non-necessarily
Hausdorff manifold $G$ with a groupoid structure over the manifold of objects
$G_{0}=M$ of $G$, given by the structure maps
$$
\xymatrix{
G\times^{st}_{M}G \ar[r]^-{\mlt} & G \ar[r]^-{\inv} &
  G \ar@<2pt>[r]^{s} \ar@<-2pt>[r]_{t} & M \ar[r]^-{\uni} & G
}
$$
which are all local diffeomorphisms.
Here $s$ and $t$ are the source and the target map, while
$\mlt$ is the partial multiplication, defined on the fibered product
$G\times^{st}_{M}G$ taken with respect to the source map on the first factor
and with respect to the target map on the second. As usual, we write $\mlt(g,h)=gh$.
We shall identify $M$ with a submanifold of $G$ via the unit map $\uni$, which
maps a point  $x\in M$ to the unit $1_{x}\in G$ at $x$.
Any arrow $g\in G$ is invertible with the inverse $\inv(g)=g^{-1}$.

If $H$ is another \'{e}tale Lie groupoid over $M$, then a smooth map $\phi:G\ra H$
is a {\em morphism of groupoids over} $M$ if $\phi|_{M}=\id$,
$s\com\phi=s$, $t\com\phi=t$ and $\phi\com\mlt=\mlt\com (\phi\times\phi)$.
We shall denote by
$$ \EtGrpd(M) $$
the category of \'{e}tale Lie groupoids over $M$ and morphisms between them.

Let $G$ be an \'{e}tale Lie groupoid over $M$.
A {\em bisection} of $G$ is an open subset $U$ of $G$ which
is both $s$-elementary and $t$-elementary. An example of a bisection is $M\subset G$.
All bisections of $G$ form a basis of topology of $G$.
Any bisection $U$ of $G$ gives a diffeomorphism
$$ \tau_{U}:t\com (s|_{U})^{-1}:s(U)\lra t(U)\;.$$
If $U$ and $V$ are two bisections of $G$, then
$U\times^{st}_{M}V$ is $\mlt$-elementary, its image
$UV=\mlt(U\times^{st}_{M}V)$
is also a bisection of $G$, and
$$ \tau_{UV}=\tau_{U}\com \tau_{V}|_{\tau_{V}^{-1}(s(U))}\;.$$
Moreover, $U^{-1}=\inv(U)$ is also a bisection of $G$ and
$$ \tau_{U^{-1}}=\tau_{U}^{-1}\;.$$

Now consider the space $\Cc(G)$ of smooth functions with compact support on $G$.
First note that because the bisections of $G$ form a basis of $G$, any
element of $\Cc(G)$ can be written as a sum of elements of the form
$$ f\com t|_{U}\;,$$
where $U$ is a bisection of $G$ and $f\in\Cc(t(U))$. Of course, the same is true if
we replace the target map with the source map, but we shall
prefer the target map.
Also, there are two natural $\Cc(M)$-module structures
on $\Cc(G)$, the one given by the source map and another given by the target map.
Indeed, both these maps are local diffeomorphisms and hence sheaves on $M$.
To make the difference transparent, we shall always write the
$\Cc(M)$-module structure on $\Cc(G)$ induced by the target map as
left action, while the
$\Cc(M)$-module structure on $\Cc(G)$ induced by the source map will be denoted as
right action.
Furthermore, the target sheaf $t:G\ra M$ on $M$ gives $\Cc(G)$ a coalgebra structure
over the left action of $\Cc(M)$, consisting of (left) $\Cc(M)$-linear maps
$$ \cm:\Cc(G)\lra \Cc(G)\tenllM \Cc(G) $$
and
$$ \cu:\Cc(G)\lra \Cc(M)\;.$$
There is another coalgebra structure induced by the source sheaf, but we will not
have a special notation for it because we will not use it here explicitly.

Now we define the {\em convolution product} $\mu$ on $\Cc(G)$ as the composition
$$
\xymatrix{
\Cc(G)\tenrlM \Cc(G) \ar[r]^-{\Omega_{s,t}}
& \Cc(G\times^{st}_{M}G) \ar[r]^-{\mlt_{+}} & \Cc(G)\;.
}
$$
Thus for any $u,v\in\Cc(G)$ we have
$$ uv=\mlt_{+}(\Omega_{s,t}(u\ten v))\in\Cc(G)\;.$$
Suppose that $U$ and $V$ are two bisections of $G$, and let
$f\in\Cc(t(U))$ and $p\in\Cc(t(V))$. Then we have
\begin{eqnarray*}
(f\com t|_{U})(p\com t|_{V})
& = & \mlt_{+}(\Omega_{s,t}(f\com t|_{U}\ten p\com t|_{V})) \\
& = & \mlt_{+}(\Omega_{s,t}(f\com\tau_{U}\com s|_{U}\ten p\com t|_{V})) \\
& = & \mlt_{+}((f\com\tau_{U})p\com s\com\pr_{1}|_{U\times^{st}_{M}V}) \\
& = & \mlt_{+}((f\com\tau_{U})p\com \tau_{U}^{-1}\com t\com\mlt|_{U\times^{st}_{M}V}) \\
& = & f(p\com\tau_{U}^{-1})\com t|_{UV}\;.
\end{eqnarray*}
Here note that $f(p\com\tau_{U}^{-1})\in\Cc(t(UV))$.
One can use this multiplication formula to easily check that the convolution product
on $\Cc(G)$ is associative and that $\Cc(M)$ is a (commutative) subalgebra of $\Cc(G)$.
Moreover, the left and the right $\Cc(M)$-action on $\Cc(G)$ coincide with
the product in $\Cc(G)$ from the left and from the right respectively.
Indeed, if $f\in\Cc(M)$ then
\begin{eqnarray*}
(f\com t|_{M})(p\com t|_{V})
& = & f(p\com\tau_{M}^{-1})\com t|_{MV}=fp\com t|_{V} \\
& = & f(p\com t|_{V})\;,
\end{eqnarray*}
while
\begin{eqnarray*}
(p\com t|_{V})(f\com t|_{M})
& = & p(f\com\tau_{V}^{-1})\com t|_{VM}=p(f\com\tau_{V}^{-1})\com\tau_{V}\com s|_{V} \\
& = & (p\com\tau_{V})f\com s|_{V}=(p\com\tau_{V}\com s|_{V})f \\
& = & (p\com t|_{V})f\;.
\end{eqnarray*}
We can therefore conclude that $\Cc(G)$ is an algebra with subalgebra $\Cc(M)$.
Note that $\Cc(G)$ has local identities in $\Cc(M)$. This algebra is called
the {\em convolution algebra} of the \'{e}tale Lie groupoid $G$, and coincide with the
Connes
definition of convolution algebra of $G$ in Hausdorff case (c.f. \cite{Con78,Con82}).

Finally, we define the antipode $S$ on $\Cc(G)$ as
$$ S=\inv_{+}:\Cc(G)\lra\Cc(G)\;.$$
If $U$ is a bisection of $G$ and $f\in\Cc(t(U))$, we have
$$ S(f\com t|_{U})=f\com t\com\inv|_{U^{-1}}=f\com s|_{U^{-1}}
   =(f\com\tau_{U})\com t|_{U^{-1}}\;.$$

One can now easily check that this gives a cocommutative coassociative
Hopf algebroid structure on $\Cc(G)$ over $\Cc(M)$.
Moreover, if $H$ is another \'{e}tale Lie groupoid over $M$ and $\phi:G\ra H$
a morphism of \'{e}tale groupoids over $M$, then
$\phi_{+}:\Cc(G)\ra\Cc(H)$ is a homomorphism of Hopf algebroids over $\Cc(M)$.
Indeed, since $\phi|_{M}=\id$ we have $\phi_{+}|_{\Cc(M)}=\id$.
We already know that $\phi_{+}$ is a homomorphism of coalgebras over $\Cc(M)$.
Furthermore, we have
\begin{eqnarray*}
\phi_{+}\com\mu
& = & \phi_{+}\com\mlt_{+}\com\Omega_{s,t}
    =\mlt_{+}\com (\phi\times\phi)_{+}\com\Omega_{s,t} \\
& = & \mlt_{+}\com\Omega_{s,t}\com (\phi_{+}\ten\phi_{+})
    =\mu\com (\phi_{+}\ten\phi_{+})\
\end{eqnarray*}
and
$$ \phi_{+}\com S=\phi_{+}\com\inv_{+}=(\phi\com\inv)_{+}=(\inv\com\phi)_{+}
   =\inv_{+}\com\phi_{+}=S\com\phi_{+}\;.$$
Therefore we proved:

\begin{prop} \label{prop18}
For any (non-necessarily Hausdorff) \'{e}tale Lie groupoid $G$
over a Hausdorff manifold $M$,
the convolution algebra $\Cc(G)$ of smooth functions with compact support on $G$
has a natural structure
of a cocommutative coassociative Hopf algebroid over $\Cc(M)$, and this gives a functor
$$ \Cc:\EtGrpd(M)\lra\HoAlgd(\Cc(M)) \;.$$
\end{prop}

We can also explicitly characterize the $S$-invariant weakly grouplike
elements of groupoid Hopf algebroids:

\begin{prop} \label{prop19}
Let $G$ be an \'{e}tale Lie groupoid over a Hausdorff manifold $M$. An element
$u\in\Cc(G)$ is $S$-invariant weakly grouplike
if and only if $u\in\Cc(U)$ for some
bisection $U$ of $G$.
\end{prop}
\begin{proof}
By Proposition \ref{prop4} it follows that there exists
a $t$-elementary $Y\subset G$ such that
$\supp(u)\subset Y$. But since we also have
$S(u)\in G^{S}_{w}(\Cc(G))$, the same proposition implies that
there exists a $t$-elementary $X\subset G$ such that
$\supp(S(u))\subset X$.
It follows that $\inv(X)$ is an
$s$-elementary neighbourhood of $\supp(u)$.
Now $U=\inv(X)\cap Y$ is a bisection with $u\in\Cc(U)$.
\end{proof}

\begin{ex}  \rm  \label{ex20}
(1) A discrete group $G$ can be viewed as an \'{e}tale Lie groupoid over
one object, $M=\{1\}$. The Hopf algebroid $\Cc(G)$ is in this case
the group Hopf algebra $\FF[G]$ over $\FF$.

(2) Let $M=\{1,\ldots,n\}$ for some $n\geq 1$, and let $G$ be the pair
groupoid $M\times M$ over $M$. In other words, for any $i,j\in M$ there is
exactly one arrow from $j$ to $i$ in $G$. The associated convolution algebra
$\Cc(M\times M)$ is the algebra of $n\times n$ matrices $\mathrm{Mat}_{n\times n}(\FF)$,
while the coalgebra structure on $\Cc(M\times M)$ is over the subalgebra of
diagonal matrices, given by
$\cm(e_{ij})=e_{ij}\ten e_{ij}$ and $\cu(e_{ij})=e_{ii}$ on the standard basis
of $\mathrm{Mat}_{n\times n}(\FF)$.
\end{ex}

\section{Spectral \'{e}tale groupoid of a Hopf algebroid}  \label{section:sgha}

Let $M$ be a Hausdorff manifold, and
let $A$ be a Hopf algebroid over $\Cc(M)$.
For any $S$-invariant weakly grouplike element $a$ of $A$,
we define an $\FF$-linear map
$$ T_{a}:\Cc(M)\lra\Cc(M) $$
by
$$ T_{a}(f)=\cu(S(fa))\;.$$
Observe that
$T_{pa}(f)=T_{a}(pf)$ for any $p\in\Cc(M)$.
Note also that we have
$$ T_{a}(f)=\cu(S(fa))=S(a)fa'=S(a')fa \;,$$
where $a'$ is any element of $A$ such that
$\cm(a)=a\ten a'$ and $\cm(S(a))=S(a')\ten S(a)$.
Indeed, we have $\cm(fa)=fa\ten a'$ and
$\cm(S(fa))=\cm(S(a)f)=S(a')\ten S(a)f=S(a')\ten S(fa)$,
so Proposition \ref{prop16} gives
$\cu(S(fa))=S(a')fa=S(fa)a'=S(a)fa'$.
If $b$ is another weakly grouplike element of $A$ with
$\cm(b)=b\ten b'$ and $\cm(S(b))=S(b')\ten S(b)$ for some $b'\in A$,
then $\cm(ba)=ba\ten b'a'$ and $\cm(S(ba))=S(b'a')\ten S(ba)$, so
$T_{a}(T_{b}(f))=S(a)S(b)fb'a'=S(ba)fb'a'=T_{ba}(f)$. Therefore we have
$$ T_{a}\com T_{b}=T_{ba}\;.$$

\begin{prop}  \label{prop21}
Let $M$ be a Hausdorff manifold, let $A$ be a Hopf
algebroid over $\Cc(M)$ and
let $a$ be an $S$-invariant weakly grouplike element of $A$
normalized on an open subset $W$ of $M$.
Then there exists a unique diffeomorphism
$$ \tau_{W,a}:V\lra W\;,$$
defined on an open subset $V$ of $M$,
such that $T_{a}(\Cc(W))=\Cc(V)$ and the restriction
$T_{a}|_{\Cc(W)}$
is given by the composition with $\tau_{W,a}$, i.e. for
any $f\in\Cc(W)$ we have
$$ T_{a}(f)=f\com\tau_{W,a}\;.$$
The $S$-invariant weakly grouplike element $S(a)$ is normalized on
$V$ and $\tau_{V,S(a)}=\tau^{-1}_{W,a}$.
The map $T_{a}|_{\Cc(W)}:\Cc(W)\ra\Cc(V)$
is an isomorphism of algebras, with the inverse
$T_{S(a)}|_{\Cc(V)}$.
\end{prop}
\begin{proof}
Take $a'\in A$ with
$\cm(a)=a\ten a'$ and $\cm(S(a))=S(a')\ten S(a)$.
First we shall prove that
$T_{a}|_{\Cc(W)}:\Cc(W)\ra\Cc(M)$ is a homomorphism of algebras.
Indeed, by Proposition \ref{prop16}
we have for any $f,p\in\Cc(W)$
\begin{eqnarray*}
T_{a}(f)T_{a}(p)
& = & S(a)fa' S(a)pa' = S(a)f\cu(a)pa' \\
& = & S(a)fpa'=T_{a}(fp)\;.
\end{eqnarray*}
Furthermore, Proposition \ref{prop16} also implies
$$ T_{S(a)}(T_{a}(f))=a'S(a)fa'S(a)=\cu(a)f\cu(a)=f \;,$$
so $T_{a}$ is injective on $\Cc(W)$.

Now choose a sequence $(f_{i})\subset\Cc(W)$
such that $f_{i}f_{i+1}=f_{i}$ for any $i\geq 1$ and
$\bigcup_{i}\supp(f_{i})=W$. Therefore
$f_{i+1}(x)=1$ for any $x\in\supp(f_{i})$,
so $\supp(f_{i+1})$ is a neighbourhood of
$\supp(f_{i})$. Define
$$ V=\bigcup_{i}\supp(T_{a}(f_{i}))\subset M\;.$$
First note that $T_{a}(\Cc(W))\subset\Cc(V)$.
Indeed,  for any $f\in\Cc(W)$ we can
choose $i$ such that $ff_{i}=f$, hence
$T_{a}(f)T_{a}(f_{i})=T_{a}(f)$ and thus
$\supp(T_{a}(f))\subset\supp(T_{a}(f_{i}))\subset V$.
Furthermore, we have
$T_{a}(f_{i})T_{a}(f_{i+1})=T_{a}(f_{i})$, and therefore
$T_{a}(f_{i+1})(x)=1$ for any $x\in\supp(T_{a}(f_{i}))$.
As before, this now implies that
$\supp(T_{a}(f_{i+1}))$ is a neighbourhood of
$\supp(T_{a}(f_{i}))$, and hence $V$ is open.
Moreover, for any $p\in\Cc(V)$ there exists and $i$
such that $T_{a}(f_{i})p=p$, and therefore
\begin{eqnarray*}
T_{a}(T_{S(a)}(p))
& = & S(a')apS(a')a \\
& = & S(a')aT_{a}(f_{i})pT_{a}(f_{i})S(a')a \\
& = & S(a')aS(a')f_{i}apS(a')f_{i}aS(a')a \\
& = & S(a')\cu(a)f_{i}apS(a')f_{i}\cu(a)a \\
& = & S(a')f_{i}apS(a')f_{i}a \\
& = & T_{a}(f_{i})pT_{a}(f_{i})=p\;.
\end{eqnarray*}
This proves that $T_{a}(\Cc(W))=\Cc(V)$, so
$$ T_{a}|_{\Cc(W)}:\Cc(W)\lra\Cc(V) $$
is an isomorphism of algebras, with the inverse $T_{S(a)}|_{\Cc(V)}$.
The open subset $V$ of $M$ is of course independent of the choice
of the sequence $(f_{i})$, because it can clearly be described
also as $V=\bigcup_{f\in\Cc(W)}\supp(T_{a}(f))$.

The isomorphism of algebras $T_{a}|_{\Cc(W)}:\Cc(W)\ra\Cc(V)$
is given by composition with
a unique diffeomorphism $\tau_{W,a}:V\ra W$.
This can be deduced from well known results of the theory
of characters on topological algebras,
however, since we are not able to give a reference, we will give
an elementary proof of this fact in Lemma \ref{lem22}
bellow.
The inverse $T_{S(a)}$ of $T_{a}$
is then of course given by composition with $\tau^{-1}_{W,a}$.

Finally, we will show that $S(a)$ is normalized on $V$.
Since $\cu(S(a))S(a')=S(a)$ by Proposition \ref{prop16},
we have for any $p\in\Cc(V)$
\begin{eqnarray*}
p\cu(S(a))\com\tau_{W,a}^{-1}
& = & T_{S(a)}(p\cu(S(a)))= ap\cu(S(a))S(a') \\
& = & apS(a)= apS(a')\cu(a) \\
& = & T_{S(a)}(p)\cu(a)=T_{S(a)}(p) \\
& = & p\com\tau_{W,a}^{-1}\;.
\end{eqnarray*}
This yields
that $\cu(S(a))$ equals $1$ on $V$. Therefore $S(a)$ is normalized on $V$ and
$\tau_{W,a}^{-1}=\tau_{V,S(a)}$.
\end{proof}

\begin{lem} \label{lem22}
Let $V$ and $W$ be Hausdorff manifolds and
$T:\Cc(W)\ra\Cc(V)$ an isomorphism of algebras. Then
$T$ is given by composition with a unique diffeomorphism
$\tau:V\ra W$.
\end{lem}
\begin{proof}
Take any point $x\in V$. Choose a sequence
$(p_{i})\subset\Cc(V)$ such that $p_{i}p_{i+1}=p_{i+1}$ and
$\bigcap_{i}\supp(p_{i})=\{x\}$. This implies that
$p_{i}(x')=1$ for any $x'\in\supp(p_{i+1})$, so
$\supp(p_{i})$ is a neighbourhood of $\supp(p_{i+1})$.
Denote $f_{i}=T^{-1}(p_{i})\neq 0$.
We have $f_{i}f_{i+1}=f_{i+1}$, so
$f_{i}(y)=1$ for any $y\in\supp(f_{i+1})$, and hence
$\supp(f_{i})$ is a neighbourhood of $\supp(f_{i+1})$.
In particular, the intersection
$$ K=\bigcap_{i}\supp(f_{i})\subset W $$
is non-empty.

We claim that in fact $K$ consists of only
one point $y$.
Indeed, assume that there are two different points
$y,y'\in K$. Choose two sequences of functions
$(w_{i}),(w'_{i})\subset\Cc(W)$ such that $w_{1}w'_{1}=0$,
$w_{i}w_{i+1}=w_{i+1}$, $w'_{i}w'_{i+1}=w'_{i+1}$, $w_{i}f_{i}=w_{i}$,
$w'_{i}f_{i}=w'_{i}$,
$\bigcap_{i}\supp(w_{i})=\{y\}$ and
$\bigcap_{i}\supp(w'_{i})=\{y'\}$. This can be done because
$\supp(f_{i})$ is a neighbourhood of $K$, for any $i\geq 1$.
Let $q_{i}=T(w_{i})\neq 0$ and $q'_{i}=T(w'_{i})\neq 0$.
We get $q_{i}q_{i+1}=q_{i+1}$ and $q_{i}p_{i}=q_{i}$, so
$$ \supp(q_{i+1})\subset \supp(q_{i})\subset\supp(p_{i})\;.$$
Analogously we get
$$ \supp(q'_{i+1})\subset \supp(q'_{i})\subset\supp(p_{i})\;.$$
This implies that
$$ \bigcap_{i}\supp(q_{i})=\bigcap_{i}\supp(q'_{i})=\{x\}\;.$$
On the other hand we have $q_{1}q'_{1}=0$ and hence also
$q_{2}q'_{2}=q_{1}q_{2}q'_{1}q'_{2}=0$.
Since $\supp(q_{1})$ and $\supp(q'_{1})$ are neighbourhoods of
$\supp(q_{2})$ and $\supp(q'_{2})$ respectively,
this implies that
$\supp(q_{2})\cap\supp(q'_{2})=\emptyset$,
a contradiction.

This proves that $K$ consists of only one point, and we denote this
point by
$$ \tau(x)\in W\;.$$
Note that this definition of $\tau(x)$ does not depend
on the choice of the sequence $(p_{i})$. Indeed, if
$(p'_{i})\subset\Cc(V)$ would be another sequence such that
$p'_{i}p'_{i+1}=p'_{i+1}$ and
$\bigcap_{i}\supp(p'_{i})=\{x\}$, we have
$p_{i}p'_{i}\neq 0$, hence $f_{i}T^{-1}(p'_{i})\neq 0$ for any
$i\geq 1$. This cannot be true unless
$\tau(x)\in\bigcap_{i}\supp(T^{-1}(p'_{i}))$.

Thus we have a well defined map
$$ \tau:V\lra W\;.$$
We have to show that $T$ is given by the
composition with this map. To this end,
take any $x\in V$, and let the sequence $(p_{i})\in\Cc(V)$ be as above,
i.e. $p_{i}p_{i+1}=p_{i+1}$ and $\bigcap_{i}\supp(p_{i})=\{x\}$,
and denote $f_{i}=T^{-1}(p_{i})$.
First observe that if $w\in\Cc(W)$
with $T(w)(x)=0$ then $w(\tau(x))=0$.
Indeed, if $w(\tau(x))\neq 0$ then $w$ has no zeros
on a neighbourhood $U\subset W$
of $\tau(x)$. We can then choose
$i\geq 1$ large enough so that $\supp(f_{i})\subset U$.
This implies that $\frac{1}{w}f_{i}\in\Cc(U)\subset\Cc(W)$, which
yields
$$ 1=p_{i}(x)=T(f_{i})(x)=T(w(\frac{1}{w}f_{i}))(x)
   = T(w)(x)T(\frac{1}{w}f_{i})(x)=0\;,$$
a contradiction. Now for any $f\in\Cc(W)$ we have
$$ T(f-T(f)(x)f_{1})(x)=T(f)(x)-T(f)(x)p_{1}(x)=0 $$
because $p_{1}(x)=1$. By the argument above
this yields
$$ (f-T(f)(x)f_{1})(\tau(x))=0\;.$$
Since $f_{1}(\tau(x))=1$, we get
$$ T(f)(x)=f(\tau(x))\;.$$

The uniqueness of such a map $\tau$ is clear by the construction.
The same argument shows that $T^{-1}$ is given by composition with
a unique map $W\ra V$, which must be the inverse of $\tau$ by uniqueness.
In particular, the map $\tau$ is a bijection.
Since the smooth maps with compact support
are preserved under the composition with $\tau$
and with its inverse as well, it follows that
$\tau$ is a diffeomorphism.
\end{proof}

\begin{cor}  \label{cor23}
Let $M$ be a Hausdorff manifold, and let $A$ be a Hopf algebroid
over $\Cc(M)$.
Suppose that $a,b\in G^{S}_{w}(A)$ are normalized on open subsets
$W$ respectively  $U$ of $M$, and let $\tau_{W,a}:V\ra W$
and $\tau_{U,b}:X\ra U$
be the corresponding diffeomorphisms.
Then for any $f\in\Cc(W)$ and any $p\in \Cc(M)$ we have
\begin{enumerate}
\item [(i)]
      $T_{a}(pf)=T_{a}(p)T_{a}(f)=(p\com\tau_{W,a})(f\com\tau_{W,a})$,
\item [(ii)]
      $T_{a}(p)|_{V}=p\com\tau_{W,a}$,
\item [(iii)]
      $fa=a(f\com\tau_{W,a})$,
\item [(iv)]
      $ba\in G^{S}_{w}(A)$ is normalized on $Z=\tau_{U,b}(W\cap X)$, and
\item [(v)]
      $\tau_{U,b}\com\tau_{W,a}|_{\tau^{-1}_{W,a}(W\cap X)}=\tau_{Z,ba}$.    
\end{enumerate}
\end{cor}
\begin{proof}
Part (i) Choose $a'\in A$ with
$\cm(a)=a\ten a'$ and $\cm(S(a))=S(a')\ten S(a)$. We have
$$ pf\com\tau_{W,a}=T_{a}(pf)=S(a)p\cu(a)fa'=S(a)pa'S(a)fa'=T_{a}(p)T_{a}(f)\;, $$
and this also yields (ii) because it holds for any $f\in\Cc(W)$.

(iii) Observe that
$$ fa=a(f\com\tau_{W,a})S(a')a=a(f\com\tau_{W,a})\cu(S(a))
   =a(f\com\tau_{W,a}) $$
because $S(a)$ is normalized on $V$.

(iv) Take $b'\in A$ with $\cm(b)=b\ten b'$ and $\cm(S(b))=S(b')\ten S(b)$.
We have
$$ \cu(ba)=baS(b'a')=baS(a')S(b')=b\cu(a)S(b')=T_{S(b)}(\cu(a)) $$
and hence (ii) implies
$$ \cu(ba)|_{Z}=\cu(a)\com\tau^{-1}_{U,b}|_{Z}\equiv 1\;.$$
This, together with $T_{a}\com T_{b}=T_{ba}$, also implies (v).
\end{proof}

Let $M$ be a Hausdorff manifold, and
let $A$ be a Hopf algebroid over $\Cc(M)$.
Since $A$ is also a coalgebra over
$\Cc(M)$, we have the associated
spectral sheaf
$\Esp(A)$ on $M$, with stalk $\Esp(A)_{y}=G(A_{y})$ for any $y\in M$.
Here we should emphasize that $A_{y}$ is obtained with respect to the left
$\Cc(M)$-action.
An \textit{arrow of} $A$ \textit{with target} $y\in M$ is
an element $g\in G(A_{y})$ for which there exists an
$S$-invariant weakly grouplike element $a$ of $A$ normalized at $y$ such that
$$ g=a_{y}\;.$$
We denote by
$$ \Gsp(A)_{y}\subset \Esp(A)_{y}$$
the set of all arrows of $A$ with target $y$.
All arrows of $A$ form a subsheaf
$$ \Gsp(A)=\mathcal{E}_{sp}^{G^{S}_{w}(A)}(A)\subset\Esp(A) $$
of the spectral sheaf $\Esp(A)$ on $M$ associated to the coalgebra $A$.
In fact, we have a basis
of open subsets of $\Gsp(A)$ of the form
$a_{W}\subset\Gsp(A)$, for $a\in G^{S}_{w}(A)$ normalized
on an open subset  $W$ of $M$.

We will now show that $\Gsp(A)$ has a natural structure of an \'{e}tale Lie
groupoid over $M$. The target map $t:\Gsp(A)\ra M$ of $\Gsp(A)$
is the projection of the sheaf $\Gsp(A)$ on $M$. Thus, if $a_{y}\in\Gsp(A)_{y}$ is
represented by an element $a\in G^{S}_{w}(A)$ normalized at $y$, we have
$$ t(a_{y})=y \;.$$
We define the source map $s:\Gsp(A)\ra M$ by
$$ s(a_{y})=\tau_{W,a}^{-1}(y) \;,$$
where $W$ is any open neighbourhood of $y$ such that
$a$ is normalized on $W$. The independence of this definition
of the choice of $a$ and $W$ clearly follows from the fact that for
another $b\in G^{S}_{w}(A)$ normalized on an open neighbourhood $U$ of $y$
there exists $p\in\Cc(M)$ which equals $1$ on an open neighbourhood
$Y\subset W\cap U$ of $y$
such that $pa=pb$, and hence
$$ \tau_{W,a}|_{\tau^{-1}_{W,a}(Y)}=\tau_{Y,a}=\tau_{Y,pa}
   =\tau_{Y,pb}=\tau_{Y,b}=\tau_{U,b}|_{\tau^{-1}_{U,b}(Y)}\;.$$
Note that $t(a_{W})=W$, while $s(a_{W})$ is precisely the domain
of the corresponding diffeomorphism $\tau_{W,a}$.

To define a partial multiplication on $\Gsp(A)$, take
any arrow $a_{y}$ with source $x$ and target $y$ and another
arrow $b_{z}$ with source $y$ and target $z$. Here
$a$ and $b$ are $S$-invariant weakly grouplike elements of $A$,
with $a$ normalized on an open neighbourhood $W$ of $y$ and $b$ normalized on an
open neighbourhood $U$ of $z$.
Define the product of these two elements by
$$ b_{z}a_{y}=(ba)_{z}\;.$$
We have to see that this definition makes sense. 
First, Corollary \ref{cor23} (iv) implies that $ba$ is normalized on
the open neighbourhood 
$Z=\tau_{U,b}(W\cap s(b_{U}))$ of $z$.
We have to show that this definition is independent on the choice of
$a$, $b$, $W$ and $U$.
First, observe that it is sufficient to show that
$(pbfa)_{z}=(ba)_{z}$ for any $f,p\in \Cc(M)$ with $f_{y}=1$ and $p_{z}=1$.
To this end, choose an open neighbourhood $W'\subset W$
such that $f$ equals $1$ on $W'$ and an open neighbourhood $U'\subset U$ of
$z$ such that $p$ equals $1$ on $U'$. Put $X=W'\cap\tau_{U,b}^{-1}(U')$,
and choose a function $\eta\in\Cc(X)$ with $\eta_{y}=1$.
In particular we have $\eta f=\eta$ and
$(\eta\com\tau^{-1}_{U,b}) p=\eta\com\tau^{-1}_{U,b}$.
Moreover, Corollary \ref{cor23} (iii) implies that
$ (\eta\com\tau^{-1}_{U,b})b=b\eta $ because $b$ is normalized on $U$ and
$\supp(\eta\com\tau^{-1}_{U,b})\subset U'\subset U$. Now
$$ (\eta\com\tau^{-1}_{U,b})pbfa=(\eta\com\tau^{-1}_{U,b})bfa=b\eta fa=b\eta a
   =(\eta\com\tau^{-1}_{U,b})ba $$
and $(\eta\com\tau^{-1}_{U,b})_{z}=1$, so $(pbfa)_{z}=(ba)_{z}$.

For any $y\in M$ there exists the unit arrow at $y$ given by
$$ 1_{y}=f_{y}\in\Gsp(A)_{y}\;,$$
where $f\in\Cc(M)\subset A$ is any function which equals $1$ on a
neighbourhood of $y$.
Any arrow $a_{y}$ with source $x$ and target $y$ has the inverse
$$ (a_{y})^{-1}=(S(a))_{x}\;.$$

\begin{theo} \label{theo24}
Let $M$ be a Hausdorff manifold and let $A$ be a Hopf algebroid over
$\Cc(M)$. Then the manifold
$$ \Gsp(A) $$
with the groupoid structure over $M$ defined above
is an \'{e}tale Lie groupoid over
$M$, called the spectral \'{e}tale groupoid of the Hopf algebroid $A$.
If $B$ is another Hopf algebroid over $\Cc(M)$ and
$\alpha:A\ra B$ a homomorphism of Hopf algebroids over $\Cc(M)$, then
$\Esp(\alpha)$ restricts to a morphism of \'{e}tale Lie groupoids
$$ \Gsp(\alpha):\Gsp(A)\lra\Gsp(B) \;,$$
and this gives a functor
$$ \Gsp:\HoAlgd(\Cc(M))\lra\EtGrpd(M)\;.$$
\end{theo}
\begin{proof}
We have the smooth structure on $\Gsp(A)$ inherited from the sheaf $\Esp(A)$,
and we just defined the groupoid structure on $\Gsp(A)$, as one can easily check.
All we need to prove is that the groupoid structure maps are local
diffeomorphisms. We already know that the target map is a local
diffeomorphism, because it is the restriction of the sheaf projection of $\Esp(A)$.
Note that the source map is also a local
diffeomorphism because on an basic open subset $a_{W}$,
given by $a\in G^{S}_{w}(A)$ normalized on an open subset $W$ of $M$,
we have
$$ s|_{a_{W}}=\tau_{W,a}^{-1}\com t|_{a_{W}} $$
and $\tau_{W,a}$ is a diffeomorphism. The same formula also gives
$$ \tau_{W,a}=\tau_{a_{W}}\;.$$
Since $\inv(a_{W})=S(a)_{s(a_{W})}$ and $t\com\inv=s$ it follows that
$\inv$ is a local diffeomorphism (in fact diffeomorphism).
If $f\in\Cc(M)$ equals $1$ on $W$ we have $\uni(W)=f_{W}$, and this
shows that $\uni$ is an open embedding.
Finally, if $b$ is another element of $G^{S}_{w}(A)$ normalized
on an open subset $U$ of $M$, we have
$\mlt(b_{U}\times_{M} a_{W})=b_{U}a_{W}=(ba)_{\tau_{U,b}(W\cap s(b_{U}))}$,
which together with $t\com\mlt=t\com\pr_{1}$ implies that $\mlt$ is a local
diffeomorphism as well.

Note that $\alpha(G^{S}_{w}(A))\subset G^{S}_{w}(B)$, which implies
that $\alpha_{y}(\Gsp(A)_{y})\subset\Gsp(B)_{y}$ for any $y\in M$.
Therefore $\Esp(\alpha)$ restricts to a morphism of sheaves on $M$
$$ \Gsp(\alpha)=\Esp(\alpha)|_{\Gsp(A)}:\Gsp(A)\lra\Gsp(B)\;.$$
We have to show that this map is in fact a morphism of groupoids.
We already know that $t\com\Gsp(\alpha)=t$, while
$\Gsp(\alpha)|_{M}=\id$ because $\alpha|_{\Cc(M)}=\id$.
If $a\in G^{S}_{w}(A)$ is normalized on an open subset $W$ of $M$,
then $\alpha(a)\in G^{S}_{w}(B)$ is also normalized on $W$,
$\Gsp(\alpha)(a_{W})=\alpha(a)_{W}$ and
$$ T_{\alpha(a)}(f)=\cu(S(f\alpha(a)))=\cu(S(\alpha(fa)))
   =\cu(\alpha(S(fa)))=\cu(S(fa))=T_{a}(f) $$
for any $f\in\Cc(M)$. This implies that
$$ \tau_{W,a}=\tau_{W,\alpha(a)} \;,$$
and hence
$$ s\com\Gsp(\alpha)|_{a_{W}}=\tau^{-1}_{W,\alpha(a)}\com t\com \Gsp(\alpha)|_{a_{W}}
   = \tau^{-1}_{W,a}\com t|_{a_{W}}=s|_{a_{W}} \;.$$
Therefore we have $s\com\Gsp(\alpha)=s$.
If $b\in G^{S}_{w}(A)$ is normalized on an open subset $U$ of $M$,
and if $z\in U$ is such that $s(b_{z})=y\in W$, then we have
\begin{eqnarray*}
\Gsp(\alpha)(b_{z}a_{y})
& = & \alpha_{z}(ba)_{z}=\alpha(ba)_{z}=(\alpha(b)\alpha(a))_{z} \\
& = & \alpha(b)_{z}\alpha(a)_{y}=\alpha_{z}(b_{z})\alpha_{y}(a_{y}) \\
& = & \Gsp(\alpha)(b_{z})\Gsp(\alpha)(a_{y})\;,
\end{eqnarray*}
so $\Gsp(\alpha)$ preserves the multiplication as well.
\end{proof}

\begin{theo}  \label{theo25}
Let $G$ be an \'{e}tale Lie groupoid over a Hausdorff manifold $M$.
Then there is a natural isomorphism of \'{e}tale Lie groupoids over $M$
$$ \Phi_{G}:G\lra\Gsp(\Cc(G)) $$
given for any $g\in G$ by
$$ \Phi_{G}(g)=(f\com \pi|_{U})_{t(g)}\;,$$
where $U$ is any bisection of $G$ with $g\in U$ and $f$
is any function in $\Cc(\pi(U))$ which satisfies $f_{t(g)}=1$.
We have $\Gsp(\Cc(G))=\Esp(\Cc(G))$, and the isomorphism $\Phi_{G}$
is precisely the isomorphism of sheaves given by Theorem \ref{theo8}.
\end{theo}
\begin{proof}
Proposition \ref{prop19} implies that $\Gsp(\Cc(G))=\Esp(\Cc(G))$,
so Theorem \ref{theo8} yields that $\Phi_{G}$ is an isomorphism
of sheaves over $M$, natural in $G$.
So the only thing we need to check is that
$\Phi_{G}$ is a morphism of groupoids.

We have $t\com\Phi_{G}=t$ because $\Phi_{G}$ is a morphism of sheaves over $M$.
Take $u\in G^{S}_{w}(\Cc(G))$ normalized on an open subset $W$ of $M$.
By Proposition \ref{prop19} we know that $u=\cu(u)\com t|_{U}$ for a bisection
$U\subset G$, and $\cu(u)$ equals $1$ on $W$.
For any $p\in\Cc(W)$ we have
$$ p\com\tau_{W,u}=T_{u}(p)=\cu(S(p\cu(u)\com t|_{U}))
   =\cu((p\com\tau_{U})\com t|_{U^{-1}})=p\com\tau_{U} $$
and therefore
$$ \tau_{W,u}=\tau_{U\cap t^{-1}(W)}\;.$$
Put $V=U\cap t^{-1}(W)$. Now we have $\Phi_{G}(V)=u_{W}$ and
$$ s\com\Phi_{G}|_{V}=\tau^{-1}_{W,u}\com t\com\Phi_{G}|_{V}
   =\tau^{-1}_{V}\com t|_{V}=s|_{V} \;.$$
This implies $s\com\Phi_{G}=s$. Observe also that
$\Phi_{G}|_{M}=\id$ because $\Phi_{G}(1_{x})=f_{x}$ for any
$f\in\Cc(M)$ with $f_{x}=1$.

Finally, we have to see that $\Phi_{G}$ preserves the multiplication.
Take any $g,h\in G$ with $s(g)=t(h)$, choose bisections $U$ and $V$
with $g\in U$ and $h\in V$, and take any
$f\in\Cc(t(U))$ with $f_{t(g)}=1$ and any $p\in\Cc(t(V))$ with $p_{t(h)}=1$.
Then
\begin{eqnarray*}
\Phi_{G}(g)\Phi_{G}(h)
& = & (f\com t|_{U})_{t(g)} (p\com t|_{V})_{t(h)}
      =((f\com t|_{U})(p\com t|_{V}))_{t(g)} \\
& = & (f(p\com\tau^{-1}_{U})\com t|_{UV})_{t(gh)}=\Phi_{G}(gh)
\end{eqnarray*}
because $UV$ is a bisection with $gh\in UV$ and
$f(p\com\tau^{-1}_{U})\in\Cc(t(UV))$
with $f(p\com\tau^{-1}_{U})_{t(g)}=f_{t(g)}p_{s(g)}=1$.
\end{proof}

\begin{theo}  \label{theo26}
Let $M$ be a Hausdorff manifold and let $A$ be a Hopf
algebroid over $\Cc(M)$.
Then there is a natural homomorphism of Hopf algebroids over $\Cc(M)$
$$ \Theta_{A}:\Cc(\Gsp(A))\lra A $$
given by
$$ \Theta_{A}(\sum_{i=1}^{k}f_{i}\com t|_{(a_{i})_{W_{i}}})
   =\sum_{i=1}^{k}f_{i}a_{i}\;, $$
where $a_{i}$ is an $S$-invariant grouplike element of $A$
normalized on an open subset
$W_{i}$ of $M$ and $f_{i}\in\Cc(W_{i})$, for any $i=1,\ldots,k$.
We have an inclusion $\Cc(\Gsp(A))\subset\Cc(\Esp(A))$
given by extension by zero, and the homomorphism $\Theta_{A}$
is the restriction $\Psi_{A}|_{\Cc(\Gsp(A))}$
of the homomorphism of coalgebras $\Psi_{A}$ given by Theorem \ref{theo9}.
\end{theo}
\begin{proof}
Denote by $\inc:\Gsp(A)\ra\Esp(A)$ the open embedding over $M$.
The map $\inc_{+}:\Cc(\Gsp(A))\ra\Cc(\Esp(A))$ is a homomorphism of
coalgebras over $\Cc(M)$ and given by extension by zero, and hence
$$ \Theta_{A}=\Psi_{A}\com\inc_{+}=\Psi^{G^{S}_{w}(A)}_{A} $$
is a homomorphism of coalgebras as well.

We have to show that $\Theta_{A}$ is a homomorphism of Hopf
algebroids over $\Cc(M)$.
First note that a function $f\in\Cc(W)$, for an open subset $W$ of $M$,
can be written as $f\com t|_{\eta_{W}}\in\Cc(\Gsp(A)$ for any
$\eta\in\Cc(M)$ which equals $1$ on $W$, and hence
$\Theta_{A}(f)=f\eta=f$. This shows that
$\Theta_{A}|_{\Cc(M)}=\id$.
Next, if $a\in G^{S}_{w}(A)$ is normalized on an open subset $W$ of $M$
and $f\in\Cc(W)$, we have
\begin{eqnarray*}
\Theta_{A}(S(f\com t|_{a_{W}}))
& = & \Theta_{A}(f\com\tau_{a_{W}}\com t|_{(a_{W})^{-1}})
      =\Theta_{A}(f\com\tau_{a_{W}}\com t|_{S(a)_{s(a_{W})}}) \\
& = & (f\com\tau_{a_{W}})S(a)=S(a(f\com\tau_{W,a})) \\
& = & S(fa)=S(\Theta_{A}(f\com t|_{a_{W}}))\;,
\end{eqnarray*}
so $\Theta_{A}$ commutes with $S$.
If $b$ is another element of $G^{S}_{w}(A)$ normalized on an open subset $U$ of $M$,
then we have for any $f\in \Cc(W)$ and any $p\in\Cc(U)$
\begin{eqnarray*}
\Theta_{A}((p\com t|_{b_{U}})(f\com t|_{a_{W}}))
& = & \Theta_{A}(p(f\com\tau^{-1}_{b_{U}})\com t|_{b_{U}a_{W}}) \\
& = & \Theta_{A}(p(f\com\tau^{-1}_{U,b})\com t|_{(ba)_{W\cap s(b_{U})}}) \\
& = & p(f\com\tau^{-1}_{U,b})ba=b((p\com\tau_{U,b})f)a=(pb)(fa) \\
& = & \Theta_{A}(p\com t|_{b_{U}})\Theta_{A}(f\com t|_{a_{W}})\;.
\end{eqnarray*}
This proves that $\Theta_{A}$ is a homomorphism of Hopf algebras.

To see that $\Theta_{A}$ is natural in $A$, take another Hopf algebroid $B$
over $\Cc(M)$ and let $\alpha:A\ra B$ be a homomorphism of Hopf algebroids over
$\Cc(M)$. The diagram
$$
\xymatrix{
\Cc(\Gsp(A)) \ar[d]_{\Gsp(\alpha)_{+}} \ar[r]^-{\inc_{+}} &
    \Cc(\Esp(A)) \ar[d]_{\Esp(\alpha)_{+}} \ar[r]^-{\Psi_{A}} & A \ar[d]^{\alpha} \\
\Cc(\Gsp(B)) \ar[r]^-{\inc_{+}} & \Cc(\Esp(B)) \ar[r]^-{\Psi_{B}} & B
}
$$
commutes because $\Psi$ is a natural transformation and $\Gsp(\alpha)$ is
the restriction of $\Esp(\alpha)$. Therefore we have
$\alpha\com\Theta_{A}=\Theta_{B}\com\Gsp(\alpha)_{+}$.
\end{proof}

\begin{theo}  \label{theo27}
For any Hausdorff manifold $M$,
the functor
$$ \Cc:\EtGrpd(M)\ra\HoAlgd(\Cc(M)) $$
is fully faithful and
left adjoint to the functor $\Gsp$.
\end{theo}
\begin{proof}
We will show that $\Phi$ of Theorem \ref{theo25} is the unit, while
$\Theta$ of Theorem \ref{theo26} is the counit of this adjunction.
Recall from Theorem \ref{theo10} that
$\Psi_{\Cc(G)}\com(\Phi_{G})_{+}=\id$ for any \'{e}ale Lie groupoid
$G$ over $M$, and that
$\Esp(\Psi_{A})\com\Phi_{\Esp(A)}=\id$ for any Hopf algebroid $A$ over $\Cc(M)$.
Since $\Gsp(\Cc(G))=\Esp(\Cc(G))$ by Theorem \ref{theo25},
we have $\Cc(\Gsp(\Cc(G)))=\Cc(\Esp(\Cc(G)))$ and
$\Theta_{\Cc(G)}=\Psi_{\Cc(G)}$, so
$$ \Theta_{\Cc(G)}\com(\Phi_{G})_{+}=\Psi_{\Cc(G)}\com(\Phi_{G})_{+}=\id \;.$$
Next, denote by
$\inc:\Gsp(A)\ra\Esp(A)$ the inclusion.
Since $\Theta_{A}=\Psi_{A}\com\inc_{+}$, it follows that
$$ \Esp(\Theta_{A})=\Esp(\Psi_{A})\com\Esp(\inc_{+})\;.$$
By definition, $\Gsp(\Theta_{A})$ is the restriction of $\Esp(\Theta_{A})$
to $\Gsp(\Cc(\Gsp(A)))$. However, Theorem \ref{theo25} implies that
$\Gsp(\Cc(\Gsp(A)))=\Esp(\Cc(\Gsp(A)))$, and therefore
$$ \inc\com\Gsp(\Theta_{A})=\Esp(\Psi_{A})\com\Esp(\inc_{+})\;.$$
This, together with the naturality of $\Phi$, now gives
\begin{eqnarray*}
\inc\com\Gsp(\Theta_{A})\com\Phi_{\Gsp(A)}
& = & \Esp(\Psi_{A})\com\Esp(\inc_{+})\com \Phi_{\Gsp(A)} \\
& = & \Esp(\Psi_{A})\com\Phi_{\Esp(A)}\com\inc=\inc\;,
\end{eqnarray*}
so $\Gsp(\Theta_{A})\com\Phi_{\Gsp(A)}=\id$.
This proves that $\Cc$ is left adjoint to $\Esp$.
Since the unit $\Phi$ of the adjunction
is a natural isomorphism by Theorem \ref{theo25}, it follows that
$\Cc$ is fully faithful.
\end{proof}

\begin{theo}  \label{theo28}
Let $M$ be a Hausdorff manifold, and let $A$ be a Hopf algebroid over
$\Cc(M)$. Then we have:
\begin{enumerate}
\item [(i)]
      $\Theta_{A}$ is an epimorphism if and only if
      $S$-invariant weakly grouplike elements of $A$ normally generate $A$.
\item [(ii)]
      $\Theta_{A}$ is a monomorphism if and only if
      $S$-invariant weakly grouplike elements of $A$
      are normally linearly independent.
\end{enumerate}
\end{theo}
\begin{proof}
This follows from Theorem \ref{theo12} with  $\Lambda=G^{S}_{w}(A)$.
\end{proof}

\begin{cor}  \label{cor29}
Let $M$ be a Hausdorff manifold and let $A$ be a Hopf algebroid over $\Cc(M)$.
If for any $y\in M$ the $(\Cs{M})_{y}$-module
$A_{y}$ is flat and if any two arrows of $A$ with target $y$ are
$(\Cs{M})_{y}$-linearly independent in $A_{y}$,
then $\Theta_{A}:\Cc(\Gsp(A))\ra A$ is a monomorphism.
\end{cor}
\begin{proof}
This follows from Theorem \ref{theo28}, Proposition \ref{prop1}
with  $L=\Gsp(A)_{y}$, and
Proposition \ref{prop11} with $\Lambda=G^{S}_{w}(A)$.
\end{proof}

\begin{theo}  \label{theo30}
Let $M$ be a Hausdorff manifold and let $A$ be a Hopf algebroid over $\Cc(M)$.
The following conditions are equivalent:
\begin{enumerate}
\item [(i)]   The Hopf algebroid $A$ is isomorphic to the groupoid Hopf
              algebroid $\Cc(G)$ of an \'{e}tale Lie groupoid $G$ over $M$.
\item [(ii)]  The homomorphism of Hopf algebroids $\Theta_{A}:\Cc(\Gsp(A))\ra A$
              is an isomorphism.
\item [(iii)] The $S$-invariant weakly grouplike elements of $A$ normally generate
              $A$ and are normally linearly independent.
\item [(iv)]  The $(\Cs{M})_{y}$-module $A_{y}$ is freely generated by
              the arrows of $A$ with target $y$, for any $y\in M$.
\end{enumerate}
Moreover, if the Hopf algebroid $A$ is isomorphic to the groupoid Hopf
algebroid $\Cc(G)$ of an \'{e}tale Lie groupoid $G$ over $M$, then
the \'{e}tale Lie groupoid $G$ is isomorphic to
the spectral \'{e}tale groupoid $\Gsp(A)$ of $A$.
\end{theo}
\begin{proof}
The equivalence between (ii), (iii) and (iv) follows directly from
Theorem \ref{theo28} and Proposition \ref{prop11} with $\Lambda=G^{S}_{w}(A)$.
The condition (ii) clearly implies (i), while
(i) implies (ii) by Theorem \ref{theo27}. The isomorphism between
$G$ and $\Gsp(A)$ is also a consequence of Theorem \ref{theo27}.
\end{proof}


\begin{thebibliography}{99}

\bibitem{CdSW}   A. Cannas da Silva and A. Weinstein,
                 {\em Geometric Models for Noncommutative Algebras}.
                 Berkeley Mathematics Lecture Notes 10,
                 American Mathematical Society, Providence (1999).

\bibitem{Con78}  A. Connes, The von Neumann algebra of a foliation.
                 {\em Mathematical Problems in Theoretical Physics,}
                 Lecture Notes in Phys. 80, Springer, New York (1978)
                 145--151.

\bibitem{Con82}  A. Connes, A survey of foliations and operator algebras.
                 Operator algebras and applications, {\em Proc. Sympos.
                 Pure Math.} 38 Part I (1982) 521--628.

\bibitem{Con94}  A. Connes, {\em Noncommutative geometry}.
                 Academic Press, San Diego (1994).

\bibitem{CraMoe} M. Crainic and I. Moerdijk, A homology theory for \'{e}tale
                 groupoids. {\em J. Reine Angew. Math.} 521 (2000) 25--46.

\bibitem{Est}    W. T. van Est, Rapport sur les $S$-atlas.
                 {\em Ast\'{e}risque} 116 (1984) 235--292.

\bibitem{Hae58}  A. Haefliger, Structures feuillet\'{e}es et cohomologie \`{a}
                 valeur dans un faisceau de groupo\"{\i}des.
                 {\em Comment. Math. Helv.} 32  (1958) 248--329.

\bibitem{Hae01}  A. Haefliger, Groupoids and foliations.
                 Groupoids in analysis, geometry, and physics,
                 {\em Contemp. Math.} 282 (2001)  83--100.

\bibitem{Lu}     J.-H. Lu, Hopf algebroids and quantum groupoids.
                 {\em International J. Math.} 7 (1996) 47--70.

\bibitem{Mal}    G. Maltsiniotis, Groupo\"{\i}des quantiques de base
                 non commutative.  {\em Comm. Algebra} 28 (2000) 3441--3501.

\bibitem{MilMoo} J. W. Milnor and J. C. Moore,
                 On the structure of Hopf algebras.
                 {\em Ann. of Math.} 81 (1965) 211--264.

\bibitem{Mrc99}  J. Mr\v{c}un, Functoriality of the bimodule associated to a
                 Hilsum-Skandalis map. {\em K-Theory} 18 (1999) 235--253.

\bibitem{Mrc01}  J. Mr\v{c}un, The Hopf algebroids of functions on \'{e}tale
                 groupoids and their principal Morita equivalence.
                 {\em J. Pure Appl. Algebra} 160 (2001) 249--262.

\bibitem{Ren}    J. Renault, {\em A groupoid approach to
                 $\mathit{C}^{\ast}$-algebras}. Lecture Notes in Math. 793,
                 Springer, New York (1980).
              
\bibitem{Tak}    M. Takeuchi, Groups of algebras over $A\otimes \overline A$.
                 {\em J. Math. Soc. Japan} 29 (1977) 459--492. 

\bibitem{Xu}     P. Xu, Quantum groupoids and deformation quantization.
                 {\em C. R. Acad. Sci. Paris} 326 (1998) 289--294.

\end{thebibliography}
\end{document}